\documentclass[11pt]{article}

\usepackage{graphicx}
\usepackage{caption}
\usepackage{subcaption}
\usepackage{float}
\usepackage{sbc-template}
\makeatletter
\renewcommand{\normalsize}{\@setfontsize\normalsize\@xipt{13.6}}
\makeatother
\usepackage[english]{babel}   
\usepackage[utf8]{inputenc}  
\usepackage{verbatim}
\usepackage{newunicodechar}
\usepackage{listings}
\usepackage{amsmath}
\usepackage{algorithm}
\usepackage{algpseudocode}
\usepackage{threeparttable}
\usepackage{acronym}
\usepackage{amssymb}
 \usepackage{booktabs} 
\usepackage{multirow}
\usepackage{graphicx}
\usepackage[export]{adjustbox}
\usepackage{subcaption}
\usepackage{makeidx}
\makeindex
\usepackage[T1]{fontenc}
\usepackage{amssymb}
\usepackage{amsmath,amssymb}
\usepackage{optidef}
\usepackage{amsmath,mathtools}
\usepackage{setspace}
\usepackage{multirow}
\usepackage{mathrsfs}
\usepackage{wrapfig}
\usepackage{tocloft}
\usepackage{calligra}
\usepackage{booktabs}
\usepackage{xcolor}
\definecolor{b}{HTML}{0B76A0}
\usepackage{mathptmx}
\usepackage{helvet}
\usepackage{courier}
\usepackage{geometry}
\usepackage{amsmath,amssymb}
\usepackage{tabularx}
 \usepackage{booktabs}
 \usepackage{siunitx}
\usepackage{caption}

\usepackage{setspace}
\usepackage{hyperref}
\usepackage{orcidlink}
\usepackage[style=apa,natbib=true,backend=biber]{biblatex}
\hypersetup{
    colorlinks=true,
    citecolor=[RGB]{40,110,150},
    linkcolor=[RGB]{40,110,150},
    urlcolor=[RGB]{40,110,150}
}
\usepackage{pdflscape}

\newcounter{subsubsubsection}[subsubsection]
\renewcommand{\thesubsubsubsection}{\thesubsubsection.\arabic{subsubsubsection}}
\newcommand{\subsubsubsection}[1]{%
  \refstepcounter{subsubsubsection}%
  \par\noindent{\bfseries \thesubsubsubsection\hspace{1em}#1}\par
}

\usepackage{colortbl} 

\usepackage{longtable}
\usepackage{array}

\newlength{\nompad}
\newlength{\nomgap}
\newlength{\nomsym}
\newlength{\nomdesc}
\newcommand{\nomhead}[1]{%
  \multicolumn{2}{|@{\hspace{\nompad}}l@{\hspace{\nompad}}|}%
    {\rule[-0.9ex]{0pt}{4.0ex}\textit{#1}}\\}
\usepackage{float}
\title{
Earth Observation Satellite Scheduling Models: Taxonomy, Classification and Open Challenges
}

\author{Margarida Caleiras\inst{1}\,\orcidlink{0009-0008-4815-7265}, Paulo Jorge Nascimento\inst{1}\,\orcidlink{0000-0003-4139-1100}\,Samuel Moniz\inst{1,*}\,\orcidlink{0000-0002-7813-4514}
}
\scriptsize{ \address{University of Coimbra, CEMMPRE, Department of Mechanical Engineering, 3030-788 Coimbra, Portugal\\
\inst{*} \ Corresponding author.\ E-mail address: samuel.moniz@dem.uc.pt (S. Moniz)}}

\begin{document} 
\onehalfspacing
\maketitle

\begin{abstract}  
Over the past three decades, a large number of optimization models have been proposed for Earth observation satellite (EOS) scheduling. The main goal of this paper is to provide an overview of these models and to relate them to the characteristics of the problems they address. To do so, we first identify the main features of EOS scheduling problems and those with the most impact on the computational complexity of the model. Based on these features, we propose a new taxonomy of satellite agility and a classification of EOS scheduling problems that allows the problem under study to be clearly identified. We then propose a classification of the optimization models into five categories, based on their time and event representations, and describe representative formulations of each category in a unified notation, together with their extensions, discussing the respective advantages and limitations. The trade-offs between the alternative representations are analyzed, showing that discrete-time formulations allow time-dependent parameters, such as transition times and image quality, to be pre-computed, which preserves the linearity of the model but makes its size depend on the adopted discretization, whereas continuous-time formulations treat these parameters as functions of the start times, which introduces nonlinearities unless simplifying assumptions are adopted. Finally, we identify the main open challenges, including the exact modeling of time-dependent transition times and the integrated scheduling of active imaging satellites.

\vspace{0.2cm}
\end{abstract}

\noindent
{\scriptsize
\setlength{\LTleft}{0pt}\setlength{\LTright}{0pt}
\renewcommand{\arraystretch}{1.1}
\begin{longtable}{|@{\hspace{\nompad}}
                   >{\raggedright\arraybackslash}p{\nomsym}@{\hspace{\nomgap}}
                   >{\raggedright\arraybackslash}p{\nomdesc}@{\hspace{\nompad}}|}

\hline
\multicolumn{2}{|@{\hspace{\nompad}}l@{\hspace{\nompad}}|}%
  {\rule[-0.9ex]{0pt}{4.2ex}\textbf{Nomenclature}}\\
\endfirsthead
\hline
\multicolumn{2}{|@{\hspace{\nompad}}l@{\hspace{\nompad}}|}%
  {\rule[-0.9ex]{0pt}{4.2ex}\textit{}}\\
\endhead
\multicolumn{2}{|@{\hspace{\nompad}}r@{\hspace{\nompad}}|}{\textit{(continued on next page)}}\\
\hline
\endfoot
\multicolumn{2}{|@{}c@{}|}{\rule[-0.9ex]{0pt}{0.9ex}}\\
\hline
\endlastfoot

\nomhead{Sets and indices}
$i, i' \in \mathbf{I}$   & satellites, $\mathbf{I} = \{1,\ldots,I\}$ \\
$j, j' \in \mathbf{J}$   & strips, $\mathbf{J} = \{1,\ldots,J\}$ \\
$k, k' \in \mathbf{K}_i$ & orbits of satellite $i$, $\mathbf{K}_i = \{1,\ldots,K_i\}$ \\
$l, l' \in \mathbf{L}$   & ground stations, $\mathbf{L} = \{1,\ldots,L\}$ \\
$n \in \mathbf{N}$       & targets, $\mathbf{N} = \{1,\ldots,N\}$ \\
$t \in \mathbf{T}$       & time points, $\mathbf{T} = \{1,\ldots,T\}$, with $T-1$ being the number of time intervals \\
$\mathbf{VW}_{ij}$       & visible time windows of satellite $i$ for strip $j$, $\mathbf{VW}_{ij} = \{vw_{ijk} \mid k \in \mathbf{K}_{ij}\}$ \\
$\mathbf{CW}_{il}$       & communication time windows of satellite $i$ for ground station $l$, $\mathbf{CW}_{il} = \{cw_{ikl} \mid k \in \mathbf{K}_{il}\}$ \\
$m \in \mathbf{OW}_{ijk}$& observation time windows of satellite $i$ for strip $j$ in orbit $k$, $\mathbf{OW}_{ijk} = \{ow_{ijkm} \mid m \in \{1,\ldots,OW_{ijk}\}\}$ \\

\nomhead{Subsets}
$\mathbf{K}_{ij}$    & orbits of satellite $i$ with respect to strip $j$, $\mathbf{K}_{ij} \subset \mathbf{K}_{i}$ \\
$\mathbf{K}_{il}$    & orbits of satellite $i$ with respect to ground station $l$, $\mathbf{K}_{il} \subset \mathbf{K}_{i}$ \\
$\mathbf{T}^o_{ijk}$ & feasible start time points for the observation of strip $j$ by satellite $i$ in orbit $k$, $\mathbf{T}^o_{ijk} \subset \mathbf{T}$ \\

\nomhead{Decision variables\textsuperscript{*}}
$\text{x}$ & binary variable indicating observation task assignment \\
$\text{y}$ & binary auxiliary variable indicating task precedence \\
$\text{z}$ & binary variable indicating downlink task assignment \\
$\text{t}$ & continuous variable indicating task start time \\
$\text{d}$ & continuous variable indicating task duration \\

\nomhead{Parameters and functions}
$\omega_j$                & weight of strip $j$ \\
$\Gamma_j$                & profit function of strip $j$ \\
$q_j$                     & minimum image quality threshold of strip $j$ \\
$(j,0),\;(j,1)$           & endpoints of strip $j$ \\
$dur_j^o,\;dur_j^d$       & observation duration and downlink duration of strip $j$ \\
$d^{\min},\;d^{\max}$     & minimum and maximum admissible observation duration \\
$vws_{ijk},\,vwe_{ijk}$   & start time and end time of visible time window $vw_{ijk}$ \\
$cws_{ikl},\,cwe_{ikl}$   & start time and end time of communication time window $cw_{ikl}$ \\
$ows_{ijkm},\,owe_{ijkm}$ & start time and end time of observation time window $ow_{ijkm}$ \\
$\mathrm{trans}$          & transition time function/parameter \\
$\mathrm{qual}$           & image quality function \\
$\bar{E}_i,\;\bar{M}_i$   & energy and memory consumption rates per unit of time of satellite $i$ \\
$E_i^{\max},\;M_i^{\max}$ & maximum onboard energy and storage capacities of satellite $i$ \\
$\delta$                  & binary parameter indicating feasibility of transitions \\

\nomhead{Superscripts}
$a$          & task type: attitude maneuver \\
$o$          & task type: observation \\
$d$          & task type: downlink \\
$o2o$        & task order: observation followed by observation \\
$o2d,\;d2o$  & task order: observation followed by downlink, and vice versa \\
$d2d$        & task order: downlink followed by downlink with respect to the satellite \\
$d2d,gs$     & task order: downlink followed by downlink with respect to the ground station 
\end{longtable}
\vspace{-1.5em}
\noindent{\footnotesize * \,The decision variables are
listed without indices, as their index structure depends on the
formulation in which they are used.}
}

\section{Introduction}
\label{sec:intro}

Space missions have attracted growing attention over the past decades, both as a means to explore and understand the Earth and the outer space, and as an increasingly important component of the global economy \citep{Ferrari2025}. In this context, Earth observation satellites (EOSs) have become indispensable tools for acquiring data across various domains, including agricultural monitoring, urban planning, disaster response, and defense applications \citep{Lemaitre2002,Chen2019,Wang2021}. Due to the increasing demand for Earth observation data, numerous satellite systems have been deployed to fulfill this need \citep{Ferrari2025}. However, since satellite resources are limited and data demands exceed their availability, efficient utilization of existing satellite resources is critical.
This has motivated extensive work on the optimization of satellite mission operations, and in particular on the EOS satellite scheduling problem, which consists of allocating over time the operational tasks of one or several satellites and ground stations so as to maximize a given objective while respecting all physical and operational constraints \citep{Bensana1999,Lemaitre2002}. Depending on the context and modeling assumptions adopted, a number of variants of the problem can be identified.

Effectively solving these problems can substantially improve operational performance, including increasing the number and value of the observation requests that are satisfied, improving the quality and timeliness of the delivered data, and making more efficient use of the limited onboard resources of the satellites \citep{Bensana1999,Wang2021,Ferrari2025}. This is especially critical for time-critical applications, such as disaster response and security operations, whose value depends on meeting strict acquisition and delivery deadlines \citep{Ferrari2025}. However, as satellite constellations continue to grow in size and capability, achieving high operational performance has become increasingly challenging \citep{Wang2021,Li2025}.

Despite three decades of research on optimization approaches for EOS scheduling, there are still a number of major challenges and important questions that remain unresolved. First, it is not clear to what extent general approaches proposed for one variant of the problem can be applied effectively to other variants. Indeed, the literature lacks a systematic analysis of how the characteristics of the problem relate to the structure of the optimization models proposed to solve it, which makes it considerably harder to transfer approaches between variants.
It also remains unclear whether existing methods are able to handle the wide range of operational features present in the literature (e.g., variable or fixed observation durations, transition times, and resource consumption), as well as different objective functions (e.g., profit maximization, makespan minimization, and energy consumption minimization).
Moreover, the strengths and limitations of the various optimization models, and the size of the instances that can realistically be solved with each of them, remain largely unknown. The goal of this paper is to provide a comprehensive review of the literature on EOS scheduling problems, with the aim of answering these questions and identifying the main challenges and future research directions in the field.

EOS scheduling problems have been described in several review papers in recent years. \citet{Wang2021} summarize the research on the agile EOS scheduling problem, proposing a classification of formulations into continuous-time and discrete-time models and discussing some problem variations regarding profit definition, multi-objective functions, and autonomous models. They also review the solution methods proposed in the literature, including exact methods, heuristics, metaheuristics, and machine learning techniques. \citet{Zhang2021} discuss the theoretical foundations and applications of the imaging scheduling problem, structuring their review around the classification of the research object, the modeling approaches, and the solution algorithms proposed in the literature. \citet{Ferrari2025} review the main variants of the satellite scheduling problem for both Earth observation and outer space missions, covering imaging, communication, and integrated scheduling, and describe the optimization approaches proposed to solve each variant. Finally, \citet{Li2025} present a review dedicated to large-scale imaging scheduling problems, describing the planning frameworks, models, and algorithms proposed for different system compositions, distinguishing between works focused on small-scale and large-scale constellations, and discussing the main challenges and future research directions regarding scalability.

Although these reviews provide a valuable understanding of the EOS
scheduling problem and its variants, they organize the literature
primarily by problem variant, application domain, or solution method.
When time representation is discussed, as in the classification
proposed by \citet{Wang2021}, the analysis is restricted to the agile
EOS imaging scheduling problem and does not cover the event
representation, that is, the way in which task sequencing is encoded.
Yet it is largely the event representation that determines the
structural trade-offs of a formulation. 
Moreover, no work adopts a unified mathematical notation that allows the formulations proposed across different studies to be directly compared.
In comparison to previous related work, we make three primary contributions. First, we propose a classification of EOS scheduling problems, as well as of the features that characterize the optimization models used to formulate them, providing a systematic analysis of the relationship between problem characteristics and model structure. Second, we review the optimization models proposed in the literature, and propose a classification based on their time and event representations. This classification shows how each representation accommodates the operational features and objective functions found in the literature, and at what cost in terms of model size and linearity, thereby clarifying which features each class of models is able to handle. Third, we describe the main variables and equations of each model category in a unified notation, and discuss the main challenges of their implementation. The unified notation allows the characteristics of the models to be compared directly and makes their structural strengths and limitations explicit. Two research questions guide our analysis:

\textbf{RQ1.} What are the main characteristics of EOS scheduling problems, and how do they impact the structure and computational complexity of the optimization models proposed to solve them?

\textbf{RQ2.} What are the main optimization models proposed to solve EOS scheduling problems, and what are their main characteristics, strengths, and limitations?

Through \textbf{RQ1}, this paper aims to characterize EOS scheduling problems by identifying the constraints that define them, as well as the main problem features that impact the structure and complexity of the resulting optimization model. We propose a new taxonomy for satellite agility and clarify the main problem concepts, such as mission scope, transition times, and observation requirements.
These features are then combined into a new proposed problem classification, which allows each problem variant addressed in the literature to be described compactly.
This not only facilitates the comparative analysis conducted throughout this review, but also provides a more rigorous definition of the main terminology used in the EOS scheduling field.
\textbf{RQ2} aims to provide a comprehensive review of the optimization models proposed in the literature. We propose a classification of these models based on their time representation (discrete or continuous) and their associated event representations, and we describe the main variables and equations of each model category. We also discuss the main structural modeling challenges, and the problem features that introduce additional modeling and computational complexity. Finally, we summarize the most pressing future research directions identified through our analysis.

The rest of the paper is organized as follows. 
Section~\ref{sec:classification_problems} describes the constraints and assumptions of EOS scheduling problems, identifies the problem features with a direct consequence on the structure of the optimization model, and proposes a classification of these problems based on a three-field notation. Section~\ref{sec:classification_models} proposes a classification of the optimization models based on their time and event representations, and Section~\ref{sec:modeling_aspects} presents the formulations proposed for each category in a unified notation, together with their extensions. Finally, Section~\ref{sec:conclusions} presents the conclusions and the main directions for future research.

\section{Classification of EOS scheduling problems}
\label{sec:classification_problems}

EOS scheduling problems exhibit a wide variety of structural characteristics that must be carefully considered when developing optimization models to solve them. This section presents a description of the constraints and assumptions considered in the EOS scheduling literature, and a systematic classification of the main problem features that have a direct impact on the structure and computational complexity of the resulting optimization model. Building on these features, it then proposes a classification of EOS scheduling problems based on a three-field notation inspired by \citet{Graham1979}.

\subsection{Problem constraints and assumptions}
\label{sec:problem_desc}
\subsubsection{Imaging scheduling}
\begin{figure}[!htbp]
\centering
\includegraphics[width=0.55\textwidth]{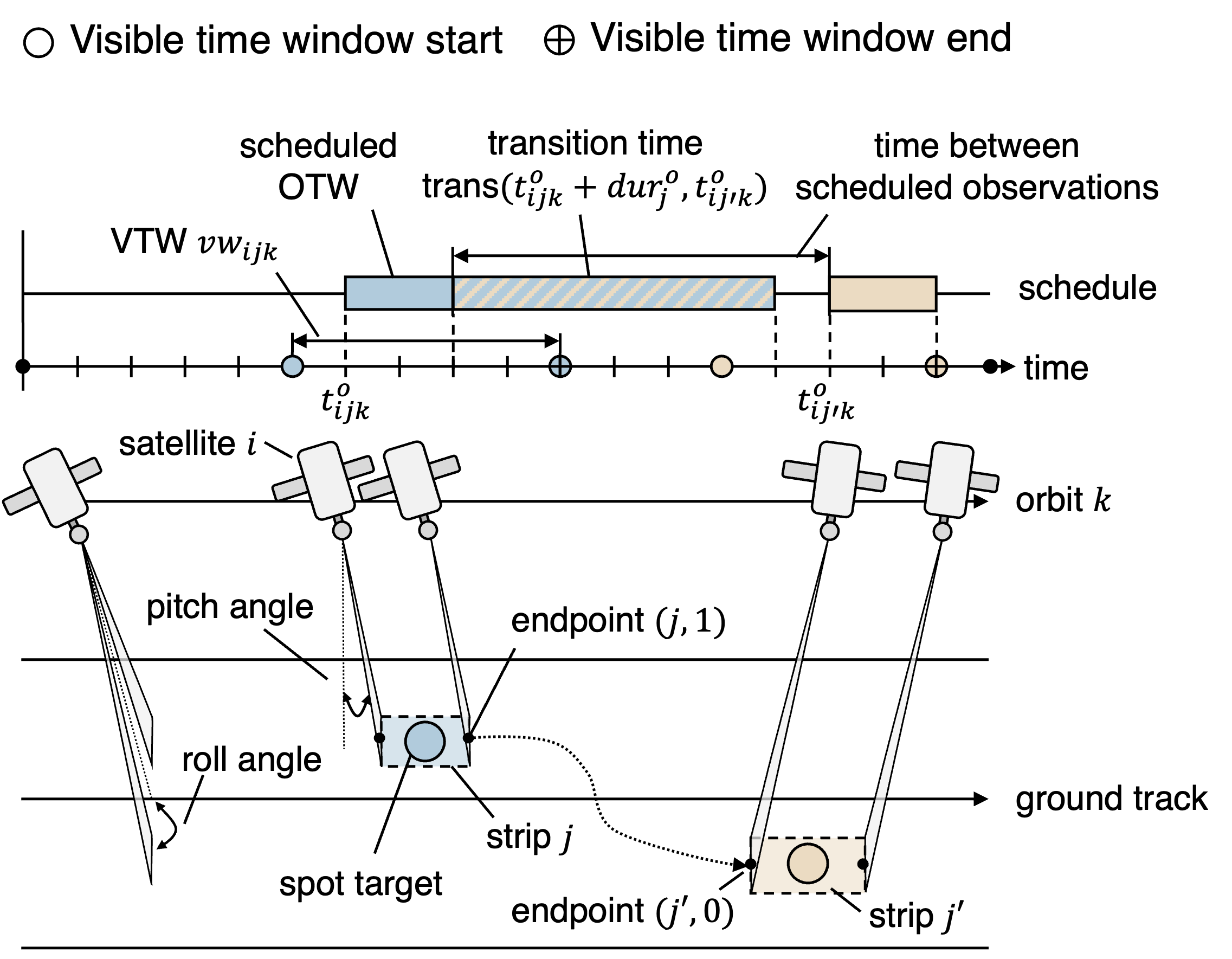}
\caption{Main concepts of imaging scheduling.}
\label{fig:imaging}
\end{figure}

\autoref{fig:imaging} illustrates the main concepts of imaging scheduling and may help the reader follow the definitions introduced throughout this section. In the example, satellite $i$ starts observing the blue strip $j$ at time $t^o_{ijk}$. After 2 units of time, the strip is observed and the satellite starts the transition to observe the yellow strip $j'$ at time $t^o_{ij'k}$. After 2 units of time the strip is observed.

Let us consider a constellation of satellites $\mathbf{I} = \{1,\ldots,I\}$ and a set of observation targets $\mathbf{N} = \{1,\ldots,N\}$.
Each target is defined by a geometrical description of the area of interest, which includes the location and the shape of the area (spot, polygon or line),
  the customer-specified due date, defined by an earliest start time and a latest deadline by which the request must be satisfied,
  and the priority of the target.
Targets are decomposed into a set of strips $\mathbf{J} = \{1,\ldots,{{J}}\}$, each strip $j \in \mathbf{J}$ with its own weight $\omega_j$, associated with the corresponding target customer's priority, as well as two endpoints $(j,0)$ and $(j,1)$.
Each satellite $i \in \mathbf{I}$ completes multiple orbits $k \in \mathbf{K}_{i}  = \{1,\ldots,{K_i}\}$ within the scheduling horizon and it can only observe a strip when it is visible from its orbit, which defines a Visible Time Window (VTW).
We denote by $\mathbf{VW}_{ij}$ the set of VTW of satellite $i$ with respect to strip $j \in \mathbf{J}$, as well as $vw_{ijk}= [vws_{ijk}, vwe_{ijk}]$ as the VTW of satellite $i$ for strip $j$ in orbit $k \in \mathbf{K}_{ij}$, where $vws_{ijk}$ and $vwe_{ijk}$ denote the start time and the end time of $vw_{ijk}$, respectively, and $\mathbf{K}_{ij}$ is the set of orbits of satellite $i$ for which strip $j$ is visible.
Each strip represents an individual observation task that must be performed continuously for a duration $dur_{j}^o$ within one of the available VTW \citep{Lemaitre2002,Bianchessi2007,Chen2019}. 
The observation duration $dur_{j}^o$ may be shorter than the full duration of the VTW, meaning the observation can start at different instants within the time window. Each $vw_{ijk}$ therefore contains a set of possible observations, called Observation Time Windows (OTW), each defined by a valid start time and the observation duration. We denote by $\mathbf{OW}_{ijk}$ the set of OTW associated with $vw_{ijk}$.

Additionally, each satellite can only perform one observation task at a time and, between two consecutively observed strips, a minimum transition time must be respected. This transition time is a function of the observation angles of the two strips, which depend on the observation timing of each strip. 
The observation angles can be defined along three axes of rotation: the roll axis, which runs along the satellite's longitudinal direction and enables lateral pointing toward a target; the pitch axis, which runs along the transverse direction and allows the sensor to move forward or backward along the ground track, which is the projection of the
orbit onto the Earth’s surface; and the yaw axis, which runs along the vertical direction and adjusts the horizontal orientation of the sensor over the Earth \citep{Lemaitre2002}. 
For each time moment $t_{ijk}$ in $vw_{ijk}$, we can obtain the satellite observation angles for each strip by calculating its relative position vector.
 Considering that the orbit period of the satellites lies between one to
several hours, the transition times for two
observations in different orbits are clearly satisfied, and we can consider only transition times for two observations scheduled in the same orbit \citep{Gabrel2002}. 
Therefore, for two consecutively observed strips $j,j' \in \mathbf{J}$ in the same orbit $k$, in VTWs $vw_{ijk}$ and $vw_{ij'k}$, at start times $t^o_{ijk},t^o_{ij'k}$ a transition time $\text{trans}(t^o_{ijk} + dur^o_{j} ,t^o_{ij'k})$ can be determined. 
Several definitions of the transition time function have been proposed in the literature, the most widely adopted being that of \citet{Liu2017}, which we define as follows. Let $roll^{(j,u)}(t)$, $pitch^{(j,u)}(t)$ and $yaw^{(j,u)}(t)$ denote the observation angles of endpoint $(j,u)$, $u \in \{0,1\}$, at time $t$. When strip $j'$ is observed immediately after strip $j$, the satellite must reorient its sensor from the last endpoint of $j$, at the end of its observation, to the first endpoint of $j'$, at the start of the following one. The total angular displacement of this transition is given by
\begin{align}
    \Delta g =& \bigl|roll^{(j,1)}(t^o_{ijk} + dur^o_j) - roll^{(j',0)}(t^o_{ij'k})\bigr| + \bigl|pitch^{(j,1)}(t^o_{ijk} + dur^o_j) - pitch^{(j',0)}(t^o_{ij'k})\bigr| \nonumber\\ &+ \bigl|yaw^{(j,1)}(t^o_{ijk} + dur^o_j) - yaw^{(j',0)}(t^o_{ij'k})\bigr|,
    \label{eq:delta_g}
\end{align}
assuming that both strips are observed from the endpoint 0 to the endpoint 1, the remaining cases being obtained by exchanging the endpoints accordingly. The transition time is then defined by the piecewise linear function
\begin{equation}
    \text{trans}(t^o_{ijk} + dur^o_j, t^o_{ij'k}) =
    \begin{cases}
        a_0,                  & \Delta g \leq 10 \\
        a_1 + \Delta g / v_1, & 10 < \Delta g \leq 30 \\
        a_2 + \Delta g / v_2, & 30 < \Delta g \leq 60 \\
        a_3 + \Delta g / v_3, & 60 < \Delta g \leq 90 \\
        a_4 + \Delta g / v_4, & \Delta g > 90
    \end{cases}
    \label{eq:trans_liu}
\end{equation}
where $\Delta g$ is expressed in degrees, $v_1, v_2, v_3$ and $v_4$ are angular velocities associated with increasing ranges of angular displacement, and $a_0, a_1, a_2, a_3$ and $a_4$ are the corresponding attitude stabilization times. The values most commonly adopted are $a_0 = 35/3$, $v_1 = 1.5^{\circ}/\text{s}$, $a_1 = 5$, $v_2 = 2^{\circ}/\text{s}$, $a_2 = 10$, $v_3 = 2.5^{\circ}/\text{s}$, $a_3 = 16$, $v_4 = 3^{\circ}/\text{s}$ and $a_4 = 22$ \citep{Liu2017,Peng2019}, for which the function is continuous and non-decreasing in $\Delta g$. Since the observation angles change continuously as the satellite moves along its orbit, $\Delta g$, and therefore the transition time, depends on the end time of the preceding observation and on the start time of the following one. An alternative definition, proposed by \citet{Gabrel1997}, computes the maneuvering time from the largest angular displacement over the three axes, and the reader is referred to that work for details.

There are also different approaches in the literature to deal with the transition times, depending on the satellite's agility and assumptions of the transition process.
Regarding the assumptions, transition times can be considered sequence-independent, meaning that the transition time between two strips is the same regardless of which two strips are observed and their sequence, and only depend on the satellite; sequence-dependent, meaning it depends on the sequence of observations; and sequence and time-dependent, meaning it does not only depend on the sequence of observations, but also on the time on which each observation occurs. If the transition time is considered sequence-independent, it is typically defined as a constant value and we have $\text{trans}(t^o_{ijk} + dur^o_{j},t^o_{ij'k}) = \text{trans}^{o2o}_i, \ \forall i \in \mathbf{I}, j,j' \in \mathbf{J}, k \in \mathbf{K}_{i}$.
 If the transition time is considered as sequence-dependent or sequence and time-dependent, it reflects the time required for the satellite to repoint its sensor from the observation angle of the last observed endpoint of strip $j$ to that of the first observed endpoint of  strip $j'$.
If the transition time is only sequence-dependent, this repointing time is treated as a fixed parameter for each pair of strips, typically precomputed at a representative instant within the VTW, and we have $\text{trans}(t^o_{ijk} + dur^o_{j}, t^o_{ij'k}) = \text{trans}^{o2o}_{ijj'k}$, or simply $\text{trans}(t^o_{ijk} + dur^o_{j}, t^o_{ij'k}) = \text{trans}^{o2o}_{ijj'}$ if the dependence on the VTW is also relaxed. If the transition time is also time-dependent, the repointing time varies with the actual scheduled observation times within the VTW, and $\text{trans}(t^o_{ijk} + dur^o_{j}, t^o_{ij'k})$ is a function of the satellite angles in $t^o_{ijk}+dur^o_{j}$ and $t^o_{ij'k}$.

Besides this, each satellite observation may need to satisfy a minimum image quality level in order to be considered successful. Similarly to transition times, image quality is a function of the observation angles that can have different definitions.
A common assumption is that observations performed closer to the nadir, which is
the angle 
between the satellite's pointing direction and the vertical 
direction towards the Earth's surface, generally yield higher image quality. The most adopted function definition based on this assumption is that of \citep{Liu2017,Peng2020}, to which the reader is referred for further details. This function represents the instant image quality of the observation at start time $t^o_{ijk}$, denoted by $\text{qual}(t^o_{ijk})$, which depends on the observation angles at the start time of the observation. 
Observation duration may also influence image quality, with longer observations increasing the camera 
exposure time and thus improving the resulting image. A definition based on this can be found in \citep{Zhou2023}. In any case, image quality is usually considered as an imposed minimum quality 
threshold $q_j$ that each observation must satisfy, and it can also be considered as an objective function, with the aim of maximizing the overall quality of the acquired images, as 
discussed further in this section.

\subsubsection{Downlink scheduling}

\begin{figure}[!htbp]
\centering
\includegraphics[width=0.55\textwidth]{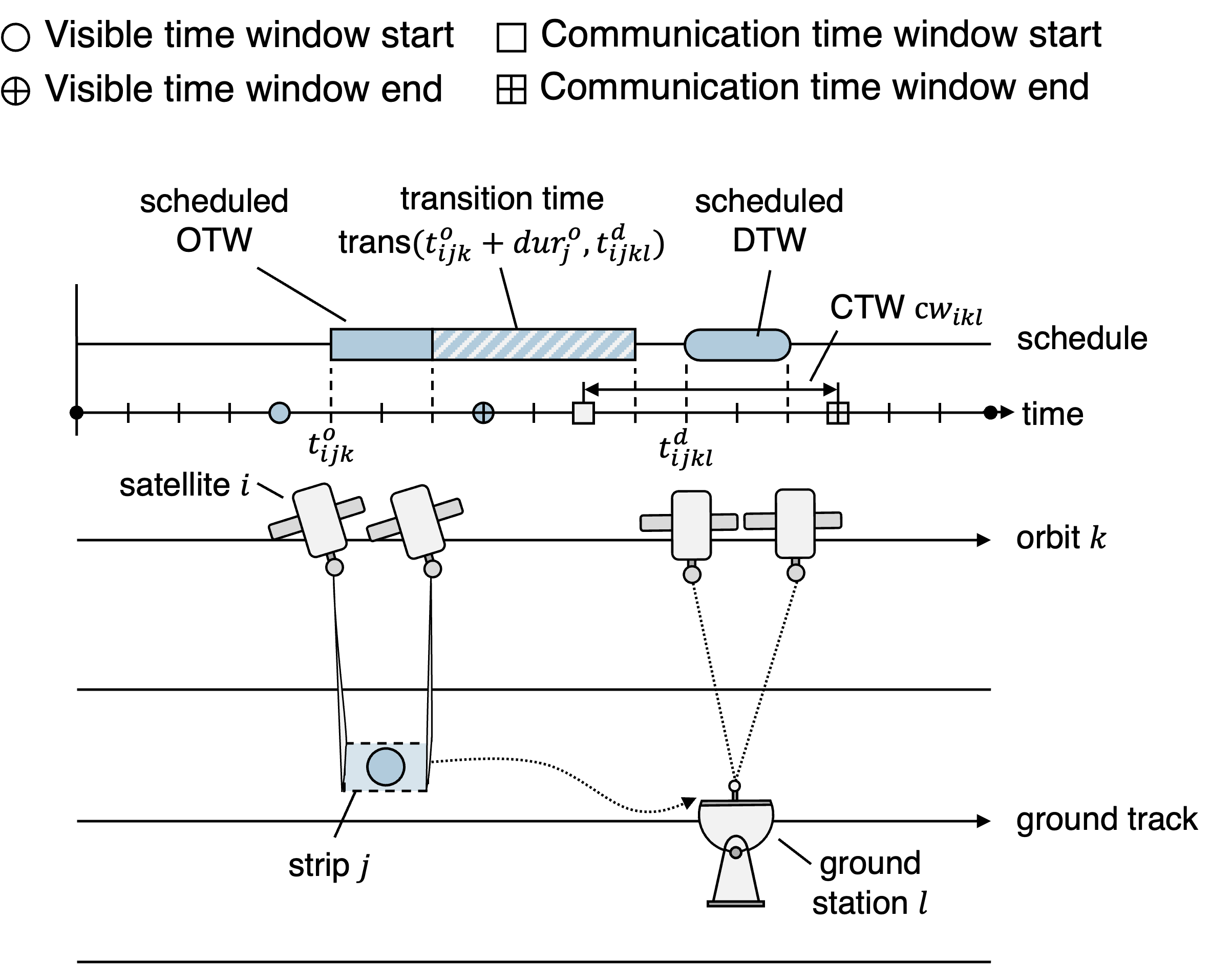}
\caption{Main concepts of downlink scheduling. }
\label{fig:downlink}
\end{figure}

\autoref{fig:downlink} illustrates the main concepts of downlink scheduling. In this example, satellite $i$ starts observing the blue strip $j$ at time $t_{ijk}^o$. After 2 units of time, the strip is observed and the satellite starts the transition to downlink the acquired data to ground station $l$ at time $t_{ijkl}^d$. After 2 units of time, the downlink is completed. 

In addition to imaging operations, each satellite must download its acquired data to a set of ground stations $\mathbf{L} = \{1,\ldots,L\}$.
Similarly to observations, a satellite can only 
communicate with a ground station when it is within its communication range, which defines a Communication Time Window (CTW). We denote by 
$\mathbf{CW}_{il}$ the set of CTW of satellite 
$i$ with respect to ground station $l \in \mathbf{L}$, and by 
$cw_{ikl} = [cws_{ikl}, cwe_{ikl}] $ the CTW of 
satellite $i$ for ground station $l$ in orbit $k \in \mathbf{K}_{il}$, 
where $cws_{ikl}$ and $cwe_{ikl}$ denote the start time and end 
time of $cw_{ikl}$, respectively, and $\mathbf{K}_{il}$ is the set of orbits of satellite $i$ during which ground station $l$ is within its communication range. The downlink corresponding to the observation of strip $j$ by satellite $i$ in orbit $k$ must be performed continuously for a duration $dur_{j}^d$ within one of the available CTW. Similarly to observations, $dur_{j}^d$ may be shorter than the full duration of a CTW, so a downlink can start at different instants within the communication window. Each $cw_{ijk}$ therefore contains a set of possible downlinks, called Downlink Time Windows (DTW), each defined by a valid start time and the downlink duration. We denote by $\mathbf{DW}_{ijk}$ the set of DTW associated with $cw_{ijk}$.

Each satellite can only perform 
one downlink operation at a time, and each ground station can communicate with
at most one satellite at a time. 
Additionally, a transition time must be considered between two consecutive operations of the same satellite and between two consecutive operations of the same ground station. That is, we must consider a transition time between two consecutive downlink operations of the same satellite to different ground stations $\text{trans}(t^d_{ijkl} + dur^d_{j},t^d_{ij'kl'}), \ j\neq j'$ and $l\neq l'$,  between a consecutive observation and downlink operation of the same satellite $\text{trans}(t^o_{ijk} + dur^o_{j},t^d_{ij'kl'})$, $\text{trans}(t^d_{ijkl} + dur^d_{j},t^o_{ij'k})$, and between two consecutive downlink operations of the same ground station $\text{trans}^{d2d,gs}$, where $t^d_{ijkl}$ denotes the start time of the downlink operation of strip $j$ to ground station $l$ by satellite $i$ in orbit $k$.
These transition times are defined analogously to the transition times between 
observations, as they depend on the angular displacement the satellite must perform 
between consecutive operations, where the ground station is treated as a single 
endpoint. The same assumptions regarding 
sequence and time-dependence apply, that is, the transitions can be treated as sequence-independent, sequence-dependent or sequence and time-dependent \citep{Zhang2022,Rocha2025}. 
We note that the transition between two consecutive downlink operations of the 
same ground station will always be treated as a constant value $\text{trans}^{d2d,gs}$, 
as its angle-based computation would depend on the ground station antenna geometry, which 
is outside the scope of this work.

\subsubsection{Resource constraints}
Each satellite is restricted by onboard resources 
limitations, such as energy and data storage capacity.
Each satellite consumes energy during imaging operations, attitude 
maneuvers and downlink operations, with observation tasks consuming $\bar{E}_i^{o}$, maneuvers consuming $\bar{E}_i^{a}$ and downlink tasks consuming $\bar{E}_i^{d}$ units of energy per unit of time, and partially recharges its onboard battery via solar 
panels collecting energy from the sun.
Although solar energy collection conditions vary due to environmental factors, the total energy collected over a single orbit remains nearly constant in practice, since orbital periods are short and illumination conditions change slowly across successive orbits. This motivates the common assumption that each satellite begins every orbit at its maximum energy capacity. This assumption explains why most papers decompose the problem by satellite orbit, treating each orbit as a separate scheduling subproblem while still accounting for the interdependencies among these subproblems. Alternatively, tasks can be indexed by VTW, and the energy level can be tracked continuously throughout the entire scheduling horizon without decomposing the problem by orbit \citep{Wang2021}. Since the majority of works in the literature that consider energy consumption adopt the orbit decomposition approach, we follow the same convention in this paper.
In this case, the energy consumption
during operations of each satellite must not exceed the maximum onboard energy capacity per orbit $E^{\max}_i$. 
Regarding data storage, each satellite is equipped with a finite 
onboard memory buffer that stores the acquired images until they 
are downloaded to a ground station. The memory level increases with 
each observation and decreases upon successful data downlink, with satellites obtaining data from observations at a rate of $\bar{M}_i^{o}$ and performing data downlinks at a rate of $\bar{M}_i^{d}$, measured in units of storage per unit of time. 
Similar to energy, the memory level must not exceed the maximum onboard storage capacity $M^{\max}_i$ \citep{Wang2011,Zhang2022}.

\begin{figure}[!htbp]
\centering
\includegraphics[width=0.55\textwidth]{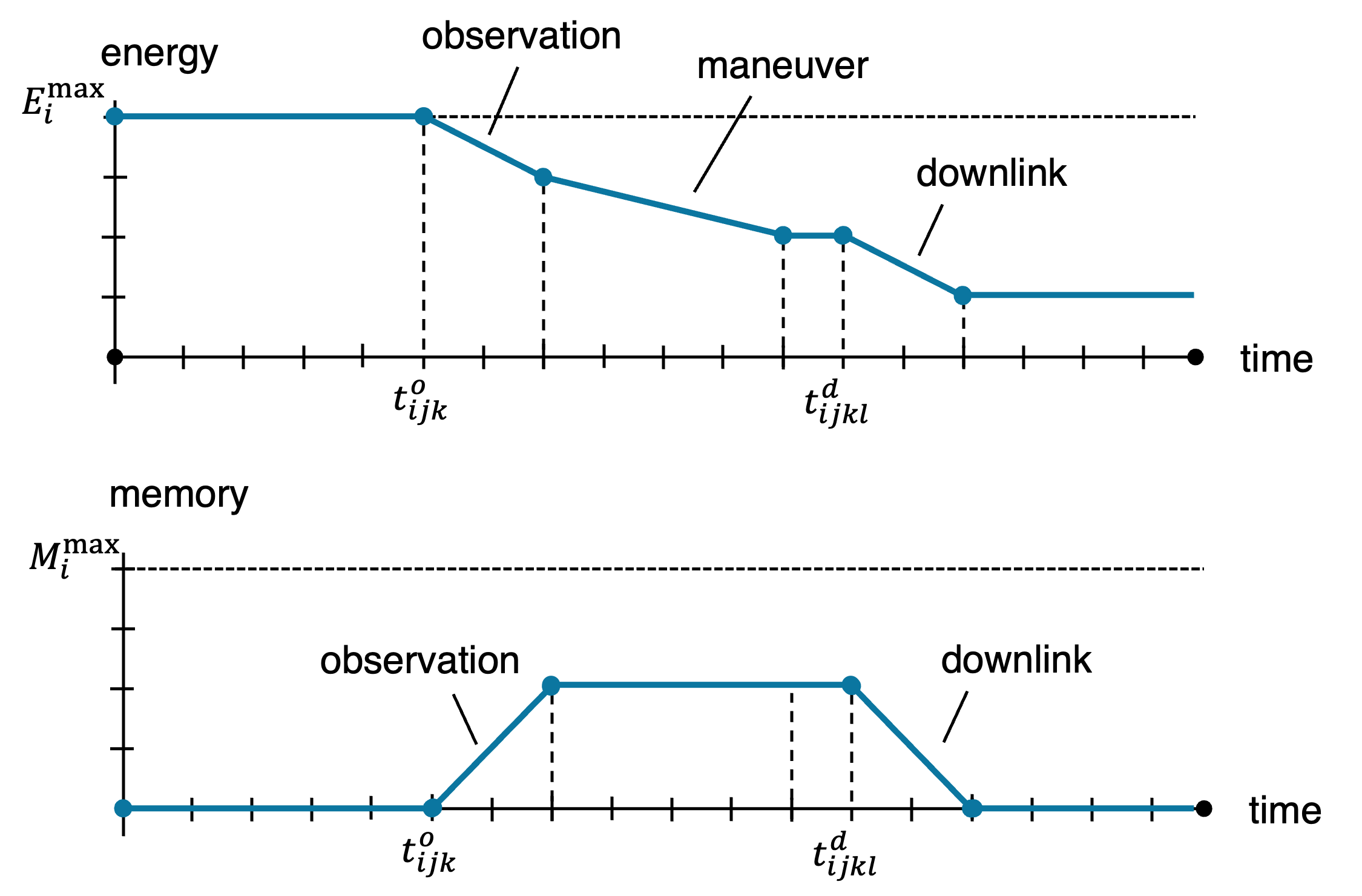}
\caption{Example of energy and memory consumption.}
\label{fig:resource}
\end{figure}

\autoref{fig:resource} illustrates the energy and memory consumption of the satellite for the example in \autoref{fig:downlink}. Here we consider that the satellite consumes
1/2 units of energy per unit of time during observations and downlinks, and 
1/4 units of energy per unit of time during maneuvers. In the example, the satellite consumes a total of
3 units of energy across the observation, the transition to the downlink, and the downlink, out of the maximum 
4 units available at the start of the first observation. Regarding memory, the satellite acquires data at a rate of $1$ unit of storage per unit of time during observations and releases it at the same rate during downlinks. In the example, the satellite acquires a total of $2$ units of data during the observation, which are stored onboard until fully transmitted during the downlink. The maximum memory capacity is $4$ units, so the memory constraint is satisfied.

\subsubsection{Objective function}
\label{sec:object_function}
The goal of the scheduling problem is to allocate observation and downlink tasks to satellites and ground stations such that all operational constraints are satisfied, while ensuring that preference is given to higher-priority requests, that the number of completed requests is as large as possible, and that the acquired data meets the required quality standards. These goals are typically formalized through an objective function, which may reflect different optimization criteria depending on the specific operational context and mission goals. 
The most common approach is to maximize a profit 
function $\Gamma_j$ associated with the scheduled observation of each strip $j$. That is, considering the binary decision variable $\text{x}_j$ that defines if strip $j$ is observed, then the objective is defined as $\max \sum_{j \in \mathbf{J}} \Gamma_j \text{x}_j$.
The profit function can take different forms, with the most simple aiming to maximize the number of satisfied requests based on the priorities of the customers, that is $\Gamma_j = \omega_j$ \citep{Lemaitre2002}. When dealing with polygonal targets, which are divided into multiple strips, it is also common to define the profit based on the percentage of the target area that is observed, penalizing partial observations. In this case, the profit function
is defined as a piecewise linear function of the observed area, with the profit associated with each strip being proportional to the area of the target that is observed \citep{Lemaitre2002, Wang2011}. 
There are also representations that introduce a more realistic characterization of 
the problem. In particular, some authors consider the so called time-dependent profit addressing the need to obtain high quality data from observations. This includes definitions of profit based on the time-dependent instant image quality function already mentioned above.
In this case, the most common profit function is defined as $\Gamma_j = \Gamma_j(t^o_{ijk}) = \omega_j \cdot 
\text{qual}(t^o_{ijk})$ \citep{Liu2017,Peng2019}.
Other common objective functions not based on the profit include makespan minimization \citep{Xiao2019}, which consists of minimizing the completion time of the last scheduled task, and energy consumption minimization \citep{Wang2014,Chang2021}. All these can be expressed as single objective or as multi-objective functions. 

\subsection{Problem features}
\label{sec:features}

\begin{figure}[!htbp]
\centering
\includegraphics[width=0.8\textwidth]{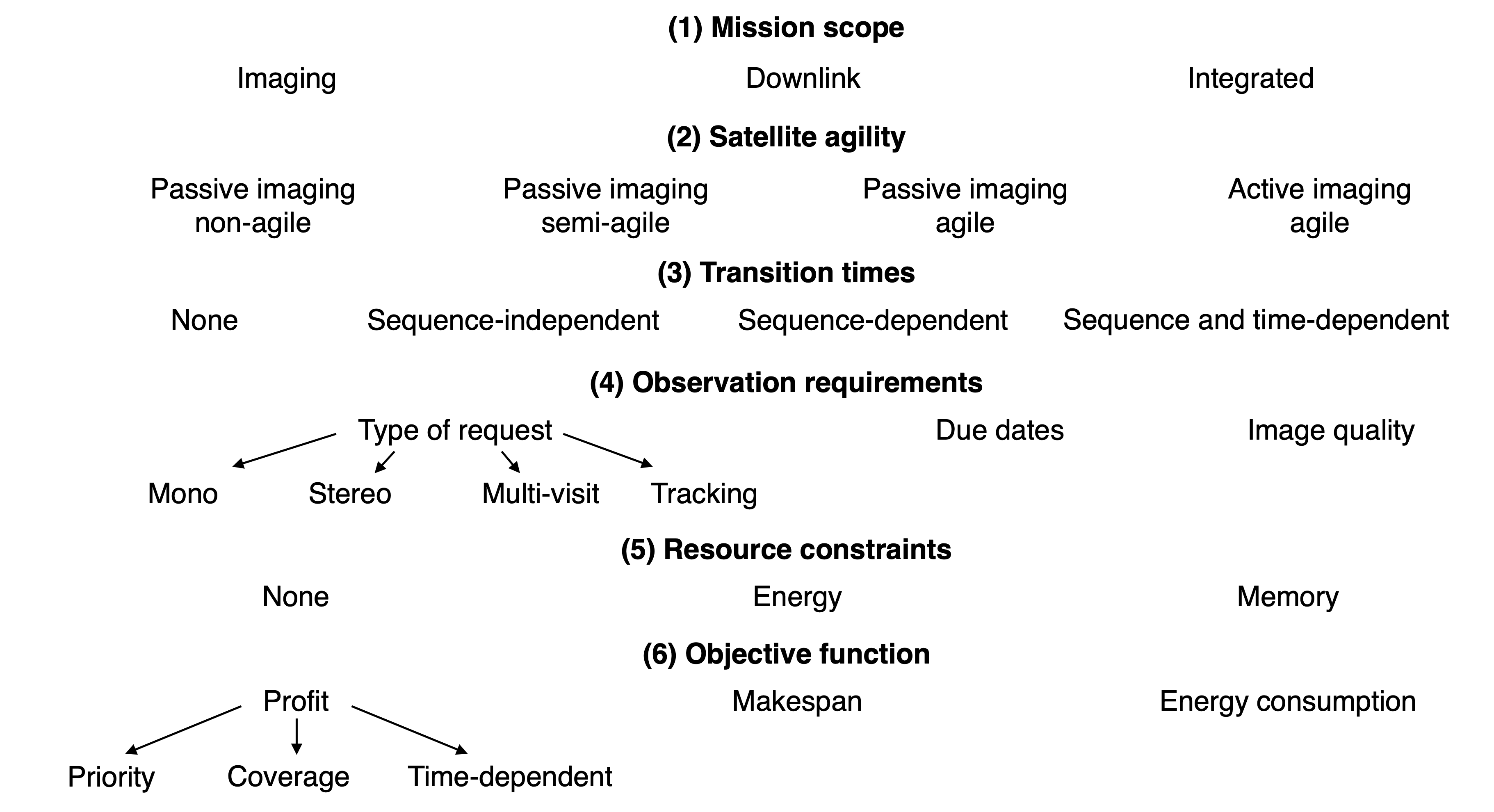}
\caption{Main features of EOS scheduling problems. }
\label{fig:features} 
\end{figure}

The following features were selected since each has a direct consequence on the structure of the optimization model, namely the types of decisions and constraints it involves and the mathematical complexity it introduces. \autoref{fig:features} summarizes the main features.

\subsubsection*{Mission scope}
First, we have the mission scope, which determines which types of tasks are explicitly scheduled and therefore which decision variables and constraints the formulation must include.
EOS scheduling problems can be classified into three main 
categories. The {imaging scheduling} deals with the selection and 
timing of observation tasks to be performed by the satellites in response 
to customer requests. The {downlink scheduling} addresses the 
allocation of communication windows between satellites and ground stations. 
Finally, the {integrated scheduling} combines both problems. This review is restricted to imaging and integrated scheduling formulations. Imaging scheduling is the most studied problem in the literature, and since integrated formulations typically adopt the same modeling structure for the observation component, it is natural to first discuss the complexities introduced by imaging scheduling and subsequently extend the analysis to the integrated imaging-downlink problem. The downlink scheduling is comparatively rare, and its most relevant modeling aspects are adequately covered within the formulations for the integrated imaging-downlink problem. For reviews focused on downlink scheduling, the reader is referred to \citep{IP2022,Ferrari2025}.

\subsubsection*{Satellite agility}
\label{sec:sat_agility}
\begin{figure}[!htbp]
\centering
\includegraphics[width=0.9\textwidth]{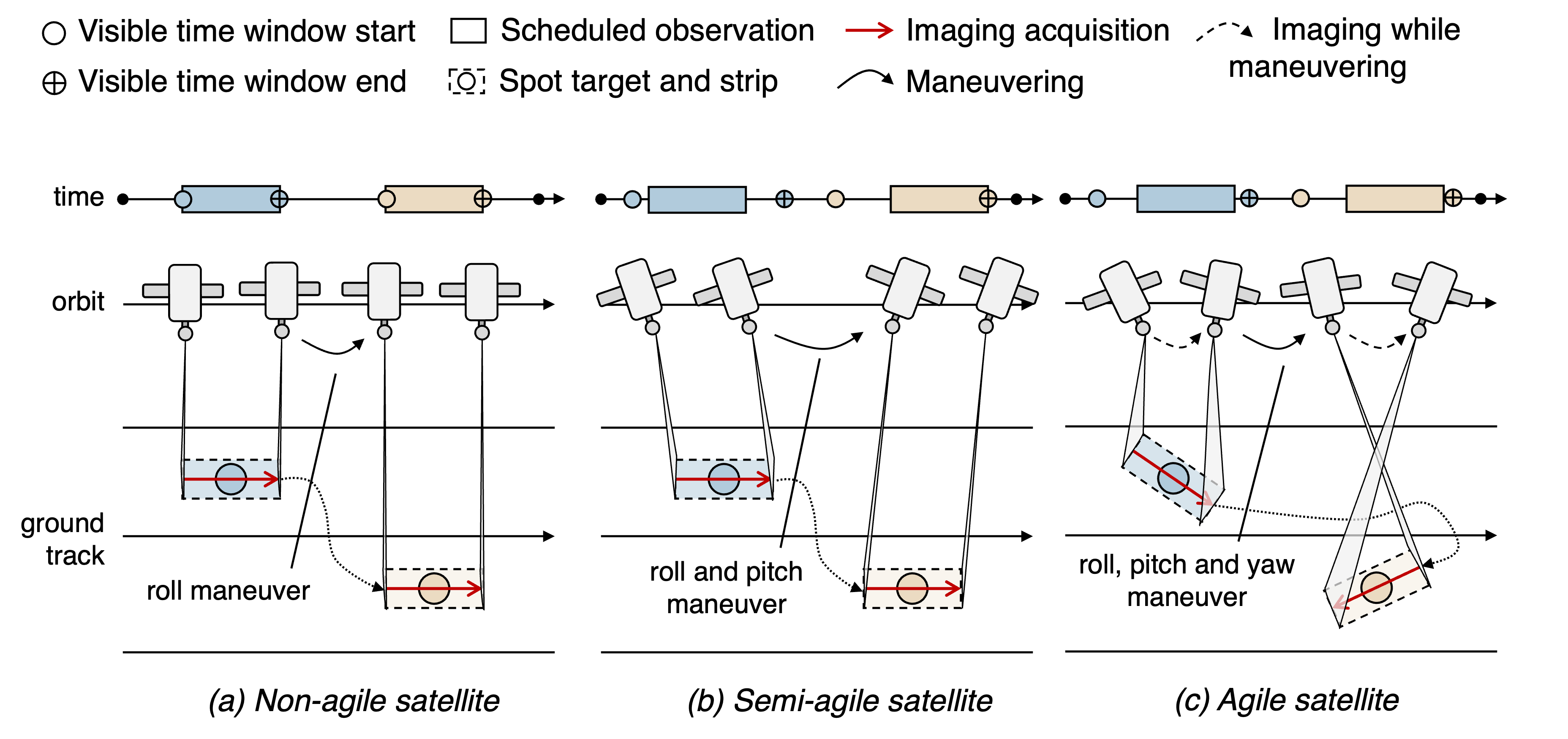}
\caption{Differences between the proposed satellite agility levels.}
\label{fig:sat_agility}
\end{figure}

A second important dimension is satellite agility, which impacts the handling of several other features of the problem, such as task start times and durations, transition times, image quality, among others. The EOS scheduling literature has traditionally classified satellites by their degree of operational flexibility, usually distinguishing between conventional, agile and super-agile satellites. However, this classification has proven ambiguous in practice. The boundaries between categories are inconsistently applied across studies, with satellites described as agile \citep{Lemaitre2002} sharing essentially the same imaging characteristics with those described as super-agile \citep{Wu2024}.
Because the conventional taxonomy provides no explicit criteria to distinguish between categories, there is no consistent basis on which to compare models and results across them, and it remains unclear which physical characteristics of each satellite drive the corresponding modeling choices.
We propose a new taxonomy grounded in the physical capabilities of the satellite and
their direct consequences for how the scheduling problem is structured. Two physical attributes are central to this characterization: axis mobility and attitude control during observations. 
Axis mobility refers to which rotation axes (roll, pitch,
yaw) the satellite can actuate, and therefore how freely it can reorient the sensor, determining both the observation opportunity windows and the set of feasible imaging
directions for each target.
Attitude control during observations refers
to whether the satellite can actively and precisely adjust the sensor while acquiring an image, as
opposed to holding a fixed attitude and relying on orbital motion. This capability
determines whether image acquisition is coupled to orbital motion or decoupled from it,
and, in the latter case, whether the imaging speed (i.e., the rate at which the
sensor's footprint sweeps across the ground) can be regulated, which in turn governs
the observation duration.
Three agility levels follow from this characterization: non-agile, semi-agile and agile satellites. \autoref{fig:sat_agility} illustrates the main operational differences between these satellites.
Non-agile satellites have mobility only along the roll axis. The image acquisition is driven by orbital motion (i.e., the sensor images targets passively as the satellite moves along its orbit, without active pointing maneuvers during the observation), meaning observations are constrained to a single direction (the orbital direction), the observation angles are fixed during observations, and the visibility windows are determined only by the satellite's ground track (see \autoref{fig:sat_agility} (a)). Semi-agile satellites have mobility along both the roll and pitch axes but still acquire images through orbital motion. Observations therefore remain constrained to a single direction and observation angles remain fixed while imaging (see \autoref{fig:sat_agility} (b)). The main difference is that pitch mobility enables look-ahead and look-back movements, extending target visibility windows and therefore introducing observation start time as a scheduling decision. Finally, agile satellites have full mobility on the three axes (roll, pitch and yaw) and are capable of stabilizing their attitude during maneuvering. This decouples image acquisition from orbital motion, enabling observations in arbitrary
directions, which requires the satellite to vary the observation angles while imaging
(see \autoref{fig:sat_agility}(c)). Furthermore, agile satellites can regulate imaging
speed by actively controlling their angular rates.

\begin{figure}[!htbp]
\centering
\includegraphics[width=0.8\textwidth]{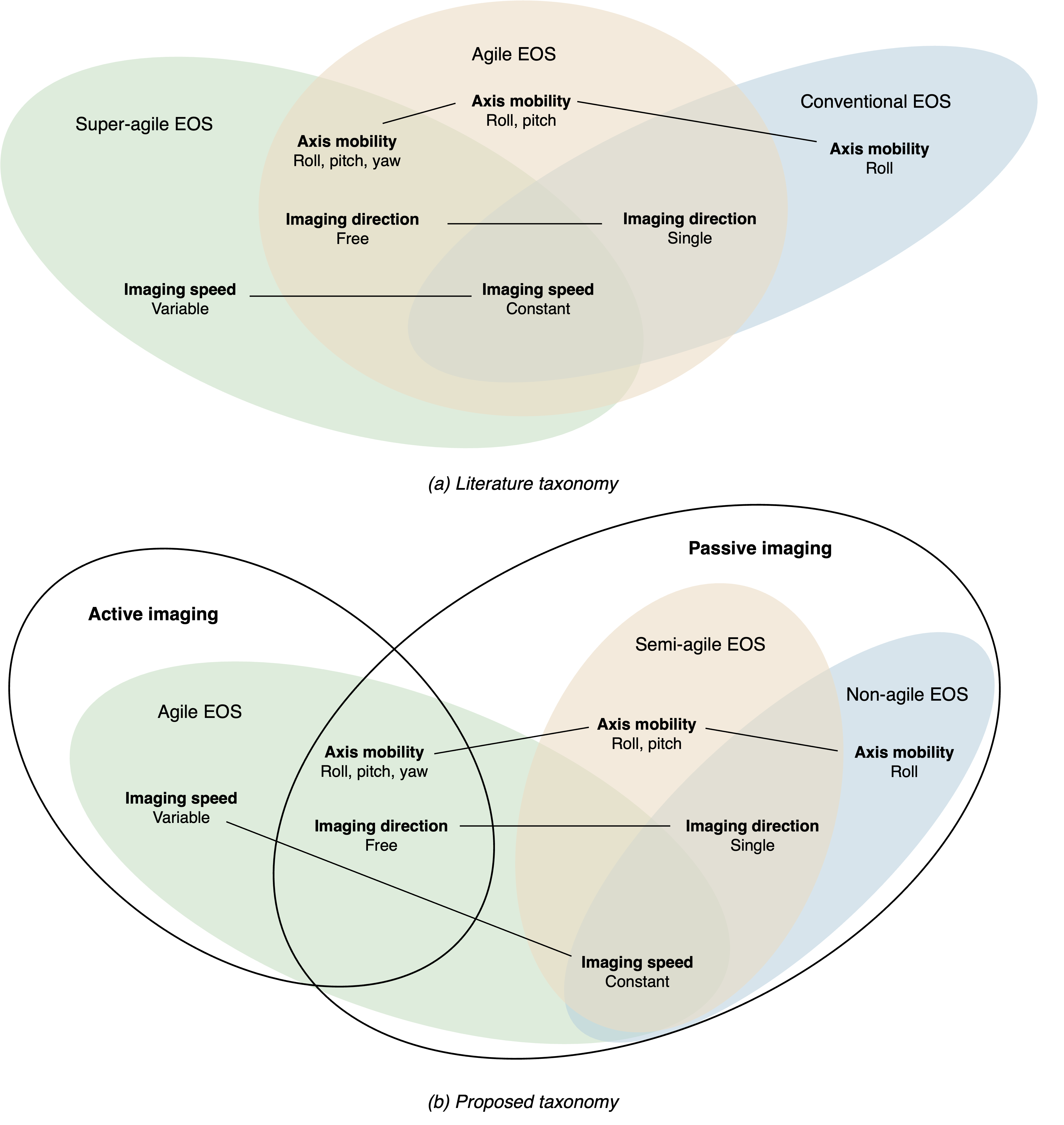}
\caption{Differences between the literature taxonomy and our proposed taxonomy. }
\label{fig:taxonomy}
\end{figure}

In the literature, only a few papers consider this regulation capability, with most authors assuming that imaging speed remains constant throughout the scheduling horizon. Since variable imaging speed introduces fundamentally different problem characteristics that directly impact both the structure of the optimization model and its computational complexity, we treat the assumption of constant imaging speed as a separate dimension of our taxonomy, which we refer to as imaging mode. We note that some authors have recognized the distinction, using terms such as active-scan \citep{Yang2018}, active imaging \citep{Chang2020} and dynamic imaging mode \citep{Lu2022a}, but it has not yet been formally established in the literature as a separate modeling dimension. 
Therefore, adopting the terminology of \citep{Chang2020}, we introduce passive imaging satellites, which cannot regulate imaging speed, and active imaging satellites, which can, introducing observation duration as an explicit optimization decision. Since active imaging requires the attitude stabilization capability that only agile satellites possess, non-agile and semi-agile satellites are necessarily passive. For agile satellites, passive and active imaging correspond to two distinct modeling assumptions. Agile satellites are referred to as passive imaging agile satellites under the assumption of constant imaging speed, and as active imaging agile satellites otherwise. 
\autoref{fig:taxonomy} summarizes the proposed taxonomy. Throughout this paper, we adopt the terminology introduced by the proposed taxonomy.

\subsubsection*{Transition times}
A third important dimension concerns transition times between consecutive tasks. As already explained, the time required to reorient the satellite between two tasks depends on the observation angles, which vary continuously as the satellite moves along its orbit. In this sense, transition times are inherently time-dependent since for the same two consecutive tasks the required reorientation time changes depending on when each task occurs. For non-agile satellites, since the set of possible start times is fully determined by the set of visibility time windows, the transition time for every possible pair of task start times can be pre-computed and treated as a fixed input parameter prior to solving. For semi-agile and agile satellites, however, task start times are a scheduling decision, making pre-computation impossible without additional simplifying assumptions.
In this context, four approaches to transition times are found in the literature: (1) transition times are not considered, (2) transition times are treated as a constant (sequence-independent), (3) transition times are treated as sequence-dependent, and (4) transition times are treated as sequence and time-dependent. 

\subsubsection*{Observation requirements}
A fourth important dimension concerns the observation requirements of each request. The first aspect is the type of request, which defines how many times and in what manner a target must be observed. A monoscopic request requires a single observation. A stereoscopic request requires two or more observations of the same target from different observation angles, enabling three-dimensional reconstruction of the scene. A multi-visit request requires repeated observations of the same target over time, as is common for disaster monitoring or time-series analysis. Finally, a tracking request involves a moving target whose position changes between visits, requiring the satellite to recompute observation opportunities at each pass and imposing additional constraints on the spacing and sequencing of observations. Formulations addressing tracking requests require significant structural modifications that fall outside the scope of this review, and the reader is referred to \citep{Chao2025,Qin2024} for relevant work in this direction.
A second aspect concerns due dates, which define the time interval within which the request must be completed. Since visibility and communication time windows are pre-computed prior to solving, due dates can be handled entirely through pre-processing by adjusting each window to fit within the due date interval. This eliminates the need to model due dates explicitly within the formulation.
Image quality is a third important requirement. The quality of an acquired image is related to the observation angles at the moment of acquisition and also on the observation duration. For passive imaging satellites, quality depends only on the observation angles at the chosen start time, which is a scheduling decision. For active imaging satellites, quality additionally depends on the observation duration, which is also a decision variable. Despite its practical relevance, image quality receives limited attention in the literature.

\subsubsection*{Resource constraints}
The fifth dimension addresses the onboard resource constraints. Two onboard resources are typically considered: energy, collected from solar panels and consumed by imaging, maneuver operations and downlinks, and memory, accumulated through observation tasks and released through downlinks. The presence of resource constraints introduces an additional feasibility dimension to the scheduling problem, since a sequence of tasks that satisfies all timing and transition requirements may still be infeasible if the resulting resource consumption exceeds the onboard capacity. 

\subsubsection*{Objective function} 
The last dimension considered in this classification concerns the objective function. Different optimization criteria can be adopted in EOS scheduling formulations, as discussed in \hyperref[sec:object_function]{Section \ref{sec:object_function}}. The choice of objective function has a direct impact on the computational performance of the model and, depending on the type of formulation adopted, may introduce additional modeling complexity, requiring extra variables and constraints for its implementation.

\subsection{Problem classification}
\label{sec:taxonomy}

The problem features described in the previous section provide the basis for a compact classification of EOS scheduling problems. Inspired by the three-field notation introduced by \citet{Graham1979} for machine scheduling, we describe each problem by a triplet $\alpha \mid \beta \mid \gamma$, in which the field $\alpha$ describes the satellite environment, the field $\beta$ describes the characteristics of the tasks to be scheduled, and the field $\gamma$ describes the optimality criteria. Each field corresponds to one or more of the features previously identified. Satellite agility, including the imaging mode, defines the satellite environment $\alpha$. Mission scope, transition times, observation requirements and resource constraints define the task characteristics $\beta$, and the objective function defines $\gamma$. Since every feature was selected for its direct consequence on the structure of the optimization model, two problems that share the same triplet can be expected to admit the same type of formulation, which makes the notation a means of identifying which formulations are applicable to a given problem variant. As in \citet{Graham1979}, we let $\circ$ denote the empty symbol, and entries equal to $\circ$ are omitted when a problem is written in this notation.

\subsubsection{Satellite environment}
The first field, $\alpha = \alpha_1 \alpha_2$, describes the satellite environment of the problem.

(1) $\alpha_1 \in \{NEOS,\; SEOS,\; PIAEOS,\; AIAEOS\}$ denotes the type of satellite, following the taxonomy of \hyperref[sec:sat_agility]{Section~\ref{sec:sat_agility}}. If $\alpha_1 = NEOS$, the satellites are non-agile; if $\alpha_1 = SEOS$, they are semi-agile; if $\alpha_1 = PIAEOS$, they are passive imaging agile satellites; and if $\alpha_1 = AIAEOS$, they are active imaging agile satellites. In practice, a constellation may combine satellites of different types, in which case the timing decisions involved would differ across satellites. Such mixed constellations are not addressed in this review, and $\alpha_1$ therefore assumes that all satellites are of the same type.

(2) $\alpha_2 \in \{\circ,\; m\}$ denotes the multiplicity of the resource. If $\alpha_2 = \circ$, the problem involves a single satellite; if $\alpha_2 = m$, it involves multiple satellites. When multiple satellites are considered, they need not be identical, and may differ, for instance, in their energy and memory capacities or in their maneuvering capabilities. These differences, however, affect mainly the parameters of the model rather than its structure, and are therefore not distinguished in this classification.

\subsubsection{Task characteristics}
The second field, $\beta \subseteq \{\beta_1, \ldots, \beta_5\}$, describes the characteristics of the tasks to be scheduled.

(1) $\beta_1 \subseteq \{O,\; D\}$ denotes the mission scope, that is, the types of task that are scheduled. If $\beta_1 = O$, only observation tasks are scheduled and the problem is an imaging scheduling problem; if $\beta_1 = \{O, D\}$, downlink tasks are also scheduled and the problem is an integrated imaging and downlink scheduling problem.

(2) $\beta_2 \in \{\circ,\; trans,\; trans_{jj'},\; trans_{jj'}(t)\}$ denotes the treatment of transition times. If $\beta_2 = \circ$, transition times are not considered; if $\beta_2 = trans$, transition times are constant (sequence-independent); if $\beta_2 = trans_{jj'}$, transition times are sequence-dependent; and if $\beta_2 = trans_{jj'}(t)$, transition times are sequence and time-dependent. When $\beta_1 = \{O, D\}$, transition times can be considered between different pairs of task types, and a different assumption may be adopted for each. In this case, $\beta_2$ takes one entry per transition type considered, each carrying the corresponding type as a superscript, following the notation $o2o$, $o2d$, $d2o$, $d2d$ and $d2d,gs$ introduced in \hyperref[sec:problem_desc]{Section~\ref{sec:problem_desc}}. A transition type that does not appear is not considered in the problem, and a single entry written without a superscript applies to all of them. For instance, $\beta_2 = \{trans^{o2o}_{jj'}(t),\; trans^{d2d,gs}\}$ denotes sequence and time-dependent transition times between consecutive observations and constant transition times between consecutive downlinks at the same ground station, with no transition considered between an observation and a downlink.

(3) $\beta_3 \in \{\circ,\; st,\; n_j\}$ denotes the type of request. If $\beta_3 = \circ$, requests are monoscopic and each strip is observed at most once; if $\beta_3 = st$, requests are stereoscopic and require two observations of the same strip, acquired from different observation angles in the same orbital direction; and if $\beta_3 = n_j$, requests are multi-visit and strip $j$ may be observed up to $n_j$ times.

(4) $\beta_4 \in \{\circ,\; q_j\}$ denotes the image quality requirements. If $\beta_4 = q_j$, each observation of strip $j$ must attain the minimum image quality threshold $q_j$; if $\beta_4 = \circ$, no quality requirement is imposed.

(5) $\beta_5 \subseteq \{E,\; M\}$ denotes the onboard resource constraints. If $E \in \beta_5$, the energy consumed in each orbit is bounded by the onboard energy capacity; if $M \in \beta_5$, the memory occupied is bounded by the onboard storage capacity. If neither is present, $\beta_5 = \circ$ and no resource constraints are considered.

\subsubsection{Optimality criteria}
The third and final field $\gamma \subseteq \{\sum \omega_j,\; \sum \Gamma^{cov}_j,\;
\sum \Gamma_j(t),\; C_{\max},\; E^{tot}\}$ corresponds to the optimality criteria chosen, following the criteria
identified in
\hyperref[sec:object_function]{Section~\ref{sec:object_function}}. If
$\gamma = \sum \omega_j$, the total priority of the scheduled strips
is maximized; if $\gamma = \sum \Gamma^{cov}_j$, the profit is defined
as a function of the observed area of each target; if
$\gamma = \sum \Gamma_j(t)$, the profit of each observation depends on
its start time through the instant image quality function; if
$\gamma = C_{\max}$, the completion time of the last scheduled task is
minimized; and if $\gamma = E^{tot}$, the total energy consumption is
minimized.
Multi-objective problems are denoted by a tuple of entries.

\subsubsection{Examples}
The following examples illustrate the use of the notation.

$SAEOS\,m \mid O,\, trans_{jj'}(t),\, \{E,\, M\} \mid \sum \omega_j$:
maximize the total profit of the observations scheduled for a constellation of
semi-agile satellites with time-dependent transition times,
subject to energy and memory limitations.

$AIAEOS \mid O,\, trans_{jj'},\, q_j \mid \sum \Gamma_j(d)$:
maximize the total time-dependent quality of the observations scheduled
for a single active imaging agile satellite, considering sequence-dependent transition times and minimum quality requirements for each observation.

$PIAEOS\,m \mid \{O,D\},\, \{trans^{o2o}_{jj'},\, trans^{d2d} \},\, M \mid
\sum \omega_j$: maximize the total profit of the observations scheduled
and downlinked by a constellation of passive imaging agile satellites with limited
memory, considering sequence-dependent transition times between consecutive observations and sequence-independent transitions between consecutive downlinks of the same satellite to differente ground stations.

\section{Classification of optimization models for EOS scheduling}
\label{sec:classification_models}
\begin{figure}[!htbp]
\centering
\includegraphics[width=0.65\textwidth]{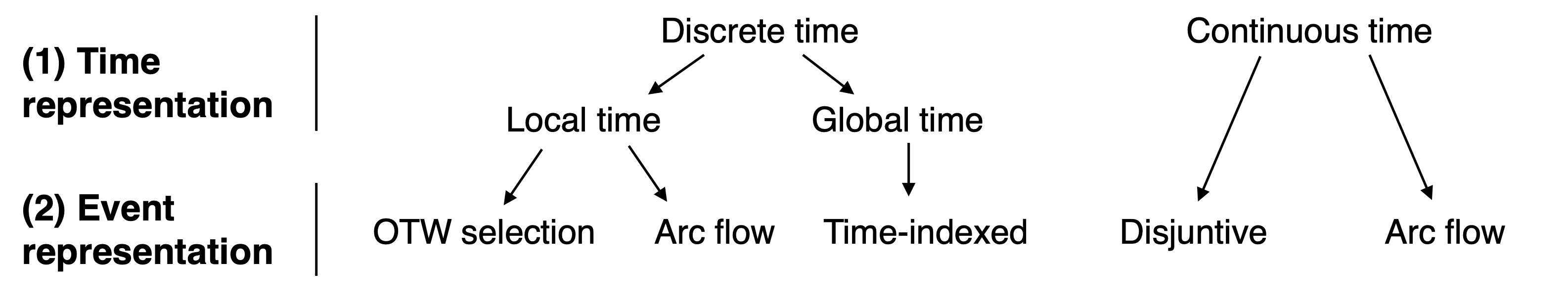}
\caption{Main dimensions of optimization models for EOS scheduling problems. }
\label{fig:event_time_tree}
\end{figure}

The classification presented here focuses on two fundamental choices that together determine the overall structure of the mathematical formulations: the time representation, which defines how timing decisions, such as task start times and durations, are handled within the model, and the event representation, which defines how the arrangement of tasks over time is encoded. These dimensions are discussed below, with particular attention to the structural trade-offs they introduce and to the ways in which they interact with the problem features described in the previous section. \autoref{fig:event_time_tree} summarizes the main categories of optimization models that arise from the different combinations of time and event representations.

\subsection{Time representation} 
    The time representation is a fundamental modeling choice in EOS scheduling, as it determines not only how the start time and duration of each task are represented within the formulation, but also how any model parameter that depends on these timing decisions is incorporated into the model.
    Depending on whether the timing
decisions are restricted to a finite set of predefined 
values or can take any value within the scheduling horizon, 
formulations can be classified as discrete or continuous.

In the context of EOS scheduling, discrete time formulations are based on: (i) local time discretization, dividing each continuous time window into a finite set of candidate tasks, each with a fixed start time and duration, and (ii) global time discretization, dividing the scheduling horizon into uniform time intervals, such that tasks must start and end at the boundaries of these intervals, which we call time points. In other words, local discretization defines a separate time grid for each visibility time window, whereas global discretization defines a single time grid shared by all tasks and satellites.
In both cases, the start time and duration of each scheduled task are implicitly determined by the chosen candidate or time interval, reducing the timing decisions to assignment decisions. 
Furthermore, since all relevant events occur at fixed and known positions, any time-dependent parameter, such as transition time, instant image quality, among others, can be pre-computed prior to solving and treated as a fixed input, which eliminates the nonlinearities that would otherwise arise from the dependence of these parameters on continuous decision variables, simplifying the mathematical structure of the formulation.
However, discrete formulations also present three important limitations. 
First, the number of candidate positions generated by the discretization depends on the chosen time interval size, which must be selected as a function of the problem data and the desired solution accuracy. Finer discretizations increase the number of candidate positions and therefore the size of the resulting model, with the extent of this increase depending on the adopted event representation.
Second, restricting timing decisions to a predefined grid may lead
to sub-optimal or even infeasible schedules.
Third, since transition times and other time-dependent parameters are computed for every pair of candidate positions prior to solving, discrete formulations incur a preprocessing cost that is absent from continuous formulations, where transition times are only evaluated for the specific task sequences that arise during the resolution process.
Despite these limitations, discrete formulations
have proven effective across a wide range of EOS operational contexts, particularly when the
number of visibility time windows is manageable and the discretization is fine enough to capture the relevant scheduling trade-offs. 

To overcome the limitations of discrete models, continuous time formulations represent 
timing decisions as explicit continuous variables, allowing tasks to 
start and end at any feasible moment within the scheduling horizon. 
This yields more flexible schedules and a model size that does not 
depend on the time interval size. However, time-dependent parameters
can no longer be pre-computed and must instead be treated as a function of 
the continuous decision variables. In particular, transition times 
between observations depend on the exact start times of the tasks, 
and since these are continuous decision variables, transition times cannot be 
evaluated prior to solving. The most common approach in the literature 
is to simplify transition times as sequence-dependent constants, 
assuming they depend only on the pair of strips involved and not on 
the exact timing of the observations. Under this simplification, the 
resulting formulation remains linear. If this simplification is not 
adopted, the timing constraints become nonlinear. Independently of this choice, the disjunctive nature of 
sequencing constraints requires big-$M$ inequalities to enforce the 
ordering of tasks.

\begin{figure}[!htbp]
\centering
\includegraphics[width=0.55\textwidth]{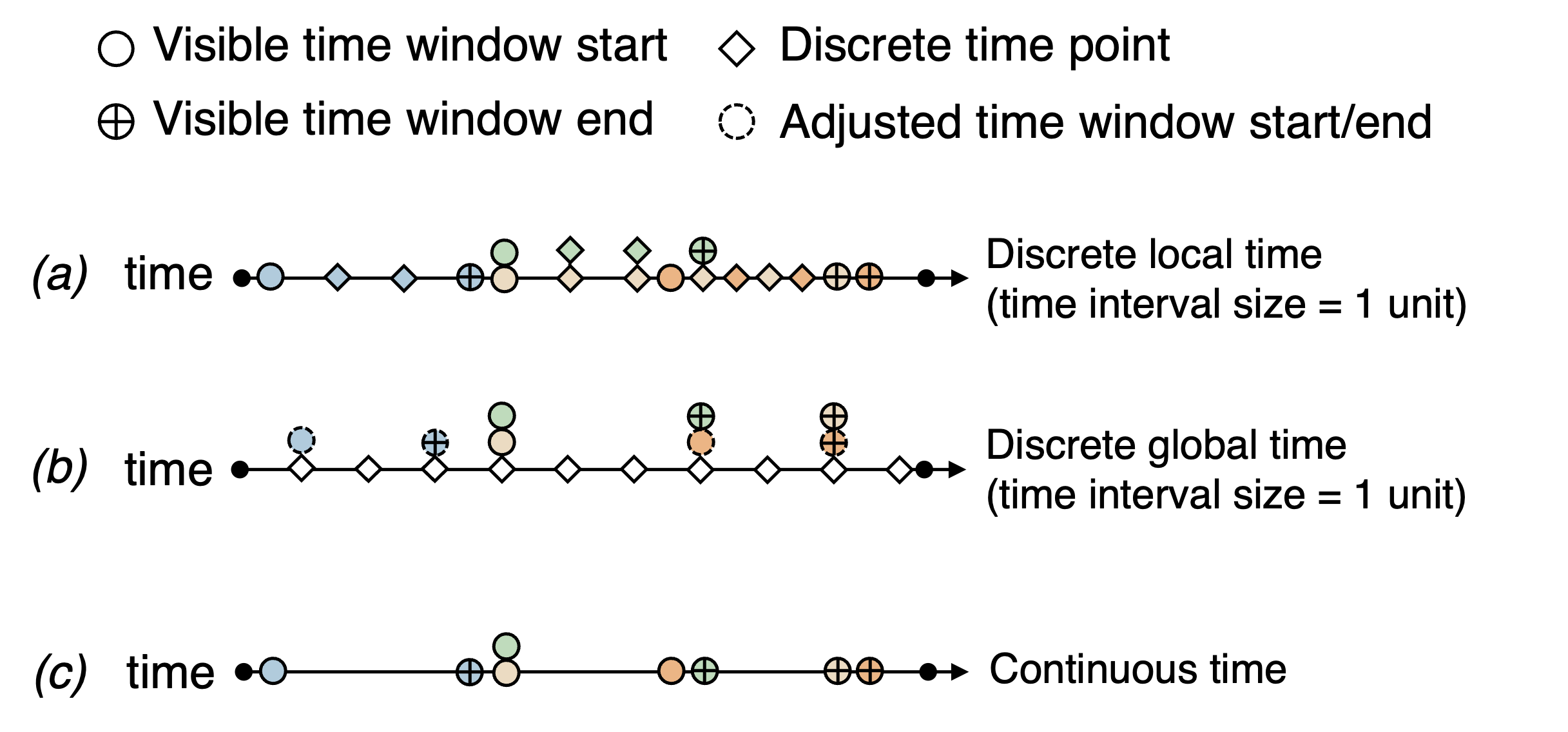}
\caption{Different time representations for EOS scheduling problems.}
\label{fig:time_representation}
\end{figure}

\autoref{fig:time_representation} illustrates the decomposition of time and the representation of the visibility windows obtained under the different time representations considered in this paper. The example involves four observation tasks, each with one associated visibility window. The three alternatives represent the problem as follows. In (a) each VTW is discretized into local time intervals, yielding 4 time points for the blue task, 4 for the green task, 6 for the yellow task, and 4 for the orange task. In (b) 10 global time points (9 time intervals) are defined, so each VTW must be anchored to this common grid. Because the window boundaries do not always coincide with the global time points, the blue and orange windows must be reduced and are divided into 2 global time points each, while the green and yellow windows are divided into 3 and 5 global time points, respectively. In (c) continuous VTWs are used, with sizes of 3, 3, 5, and 3 units of time for the blue, green, yellow, and orange tasks, respectively.

\subsection{Event representation} The event representation determines how the scheduled tasks are arranged over time within the mathematical model with the main purpose of guaranteeing that the operational constraints are satisfied and that the maximum capacity of the shared resources is never exceeded. We defined five different types of event representations that have been most used in the considered literature to develop a variety of mathematical formulations for EOS scheduling problems.

For discrete time formulations, three types of event representations are identified, which we refer to as the OTW selection, discrete arc flow and time indexed representations.
The first two are both based on local time discretization and differ in how task sequencing is represented. In the OTW selection representation, each candidate observation is treated as an independent task, and scheduling reduces to selecting a feasible subset of non-conflicting candidates. Whether two candidates can be consecutively scheduled is determined prior to solving by checking whether the transition time between them is satisfied, and this information is encoded as a binary parameter that enters the sequencing constraints directly. As a result, sequencing constraints reduce to simple linear inequalities over the binary selection variables, with no big-$M$ relaxations required. In the arc flow representation, the candidate tasks for each satellite and orbit are organized as nodes of a directed graph, and an arc between two nodes is included only if the corresponding transition is feasible. Scheduling then corresponds to finding a path through this graph, and the arc variable itself identifies the direct successor of each selected observation. The key structural difference relative to OTW selection is that infeasible transitions are not enforced through constraints but are absent from the model by construction, which eliminates sequencing constraints entirely. 
 The third discrete representation is the time indexed formulation, based on global time discretization. Unlike the previous two, this representation does not encode transition feasibility structurally, and sequencing constraints require explicit ordering variables together with big-$M$ relaxations to enforce the ordering of tasks. These ordering variables can be defined under immediate or general precedence, as discussed further ahead. Its main advantage over local discretization approaches is that the uniform time grid defines a common temporal axis shared across all tasks and satellites, which makes it easier to coordinate tasks of different types and simplifies the modeling of certain extensions, as discussed in the following sections.

For continuous time formulations, two types of event representations are found, which we refer to as the disjunctive and continuous arc flow representations. In the disjuntive representation, scheduling decisions are expressed through two separate sets of variables: binary allocation variables that determine which tasks are selected for observation, and binary ordering variables that encode the relative sequence of tasks assigned to the same satellite, which can again be defined under immediate or general precedence. In the arc flow representation, observations are modeled as nodes of a directed graph and the scheduling decisions are expressed as binary arc variables that encode immediate succession between pairs of tasks, with continuous variables representing their start times. This representation shares the directed-graph structure of the discrete arc flow representation, but since start times are decision variables, transition feasibility can no longer be checked prior to solving and must instead be enforced through timing constraints with big-$M$ parameters.

The disjuntive representation thus shares with the time indexed representation the use of explicit ordering variables to enforce transition feasibility, and in both cases these variables can be defined under immediate or general precedence.
Under immediate precedence, an ordering variable takes the value one only when one task is the direct successor of the other, which requires flow conservation constraints to ensure that each selected task has exactly one immediate predecessor and one immediate successor. 
Under general precedence, an ordering variable takes the value one whenever one task precedes another, not necessarily immediately, which avoids flow conservation constraints but requires constraints ensuring that, for every pair of tasks assigned to the same satellite, exactly one of the two possible orderings is selected.
We further note that general precedence enforces the transition time between 
every ordered pair of tasks assigned to the same resource, rather than only 
between consecutive ones. This is equivalent to the physical requirement only if 
the transition time function satisfies a triangle-type condition, namely that the 
transition between any two tasks never exceeds the transition through an 
intermediate task plus the duration of that task. The condition holds for the 
transition time functions most commonly adopted in the EOS scheduling literature, 
since the angular displacement measures on which they are based satisfy the 
triangle inequality and the associated timing functions are increasing and 
subadditive. When it does not hold, general precedence restricts the feasible 
region and may cut off optimal schedules, whereas immediate precedence remains 
exact in all cases. The choice between the two variants therefore affects not 
only the number and the type of constraints in the model, but also the conditions 
under which the formulation is exact.

\begin{figure}[!htbp]
\centering
\includegraphics[width=1\textwidth]{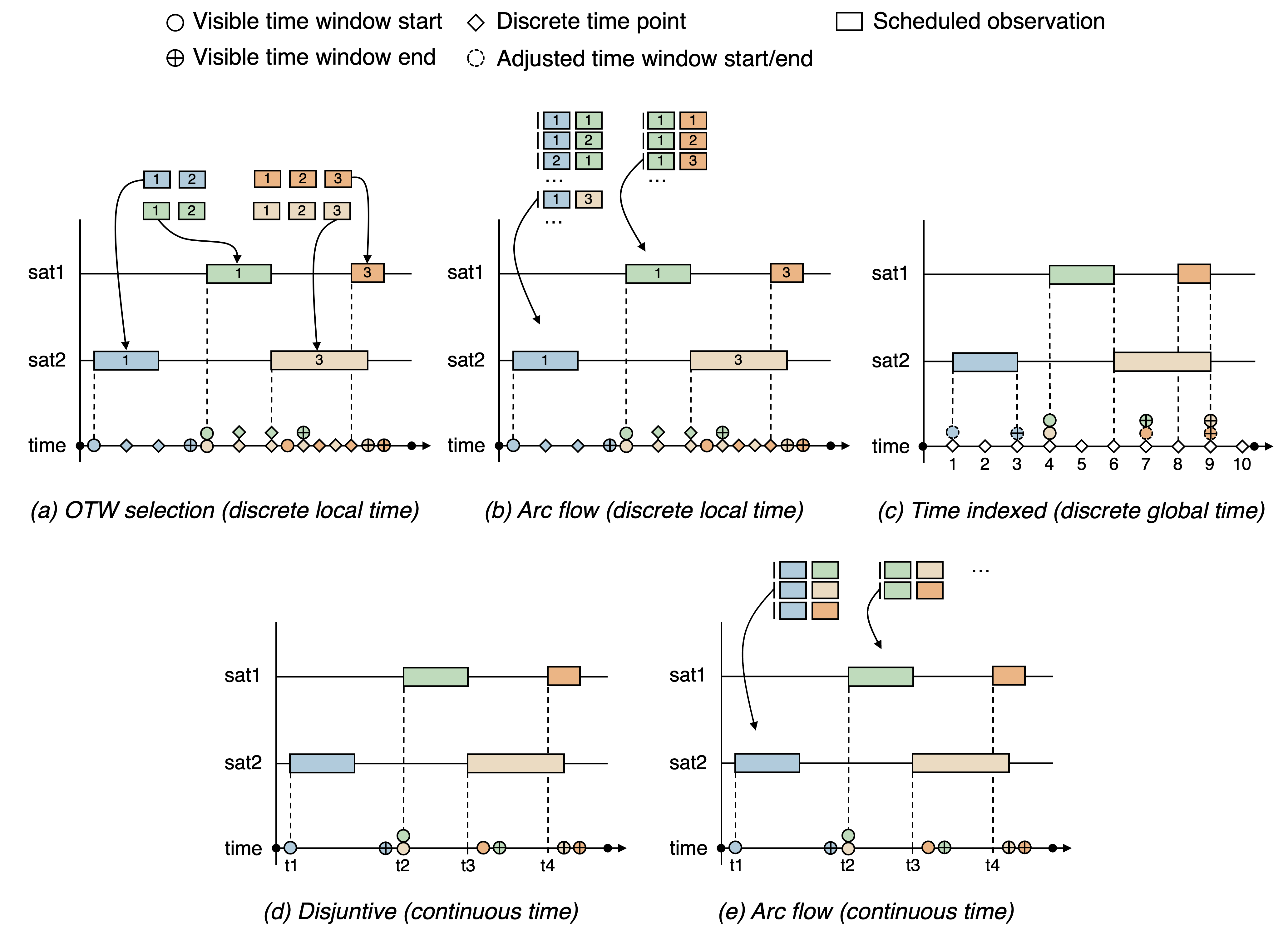}
\caption{Different event representations for EOS scheduling problems.}
\label{fig:event_time}
\end{figure}

\autoref{fig:event_time} shows a schematic representation of identical imaging schedules obtained using the alternative event representations considered in this paper. The example, based on the same problem defined in \autoref{fig:time_representation}, 
involves four observation tasks allocated to two satellites (sat1 and sat2), 
with processing durations of 2, 2, 1, and 3 units of time, respectively, each 
associated with a visible time window. For simplicity, both satellites are 
assumed to share the same visible time windows, and all transition times are 
treated as feasible. To represent the solution, the different alternatives require: (a) 4 scheduled observations selected from 10 candidates; (b) 4 scheduled observations and the 2 arcs connecting the consecutive tasks assigned to each satellite, selected from 24 candidate arcs connecting the 10 candidate tasks; (c) 4 observation start times selected from 10 candidate global time points; (d) 4 continuous observation start times; and (e) 4 continuous observation start times and the 2 arcs connecting the consecutive tasks assigned to each satellite, selected from 7 candidate arcs.

\section{Modeling aspects of alternative approaches}
\label{sec:modeling_aspects}
Having introduced a classification of EOS scheduling problems and optimization models, we present in this section the specific formulations proposed in the literature for each of the identified categories in \hyperref[sec:classification_models]{Section \ref{sec:classification_models}}. Some formulations were rewritten from their original version in order to adopt a consistent nomenclature and model structure throughout this review. 

To facilitate the comparison of formulations, we introduce a reference nomenclature for the main constraints and assumptions described above and that will be considered to model the EOS scheduling problem. \autoref{tab:constraints} shows the constraints associated with imaging scheduling (I), downlink scheduling (D), integrated scheduling (ID), and resource constraints (R) as well as the main problem assumptions (A).
The constraints can be classified into four categories that reflect their role in the mathematical formulation. Allocation constraints determine which tasks are assigned to which resources, sequencing constraints enforce the ordering of tasks on shared resources, timing constraints define when individual tasks can occur, and resource constraints limit the consumption of onboard resources throughout the mission. 

\begin{table}[!htbp]
    \small
\centering
\caption{Assumptions and constraints of EOS scheduling problems.}
\label{tab:constraints}
\begin{tabularx}{\textwidth}{@{}lX@{}}
\toprule
{Code} & {Assumptions or constraints} \\
\midrule
I1  & Each satellite can only perform one observation task at a time \\
I2  & Each strip is observed at most once (monoscopic request) \\
I3  & Each observation task must be scheduled within a visibility time window \\
I4  & A transition time is required between two consecutive observation tasks assigned to the same satellite \\
I5  & Each observation task must satisfy a minimum image quality level \\
D1  & Each satellite can only perform one downlink task at a time \\
D2  & Each ground station can only perform one downlink task at a time \\
D3  & Each downlink task can only be scheduled within a communication time window \\
D4  & A transition time is required between two consecutive downlink tasks assigned to the same satellite \\
D5  & A transition time is required between two consecutive downlink tasks assigned to the same ground station \\
ID1 & Each observed strip must also be downlinked \\
ID2 & The downlink of a strip must be scheduled after its observation has been completed \\
ID3 & A transition time is required between a consecutive observation and downlink tasks assigned to the same satellite \\
R1  & Each satellite has limited onboard energy capacity \\
R2  & Each satellite has limited onboard memory capacity \\
A1a  & The transition time is sequence-independent \\
A1b  & The transition time is sequence-dependent \\
A1c  & The transition time is sequence and time-dependent \\
A2  & The imaging speed of each satellite is assumed to remain constant throughout the scheduling horizon \\
\bottomrule
\end{tabularx}
\end{table}

We note that the classification introduced in \hyperref[sec:taxonomy]{Section~\ref{sec:taxonomy}} and the codes of \autoref{tab:constraints} describe the problem at two different levels. The classification states which features a problem includes, whereas the codes identify the constraints and assumptions through which a formulation represents them. Constraints I1 and I3 apply to every problem, and some features determine their constraints directly. Monoscopic requests, $\beta_4 = \circ$, correspond to constraint I2, a minimum image quality, $\beta_3 = q_j$, requires constraint I5, and each element of $\beta_5$ requires constraint R1 or R2. Each transition term in $\beta_2$ requires the corresponding constraint among I4, D4, D5 and ID3, under the assumption among A1a to A1c that matches its dependence, so that a single problem may combine several of these assumptions. Other features leave more freedom to the formulation. Integrated scheduling, $\beta_1 = \{O,D\}$, introduces constraints D1 to D3 and ID2, but not necessarily ID1, since some formulations penalize the data that is not downlinked instead of requiring every observed strip to be downlinked \citep{Rocha2025}. Finally, assumption A2 holds for all satellite types except active imaging agile satellites, $\alpha_1 = AIAEOS$.
Each formulation type is presented in terms of a general model and a set of extensions. The general model is defined as a base formulation for the multi-satellite imaging scheduling problem considering passive imaging semi-agile satellites, which implies fixed observation durations and a single imaging direction per strip. In this model, we will consider constraints \hyperref[tab:constraints]{I1}--\hyperref[tab:constraints]{I4}, \hyperref[tab:constraints]{R1} and \hyperref[tab:constraints]{R2}, expressed as allocation constraints, timing constraints, sequencing constraints and resource constraints. 
The optimality criterion is the maximization of the total profit of the
scheduled observations, which is the most common objective in the literature. We note that we
do not restate it for each formulation, since it is a direct adaptation
of the binary observation assignment variables of each model. Using the classification introduced in \hyperref[sec:taxonomy]{Section~\ref{sec:taxonomy}}, the general model corresponds to the problem $SEOS\,m \mid O,\, trans_{jj'}(t),\, \{E,\, M\} \mid \sum \omega_j$.

We refer to extensions as additional problem dimensions considered in the literature, namely objective function variations, image quality requirements, active imaging satellite characteristics, and integrated downlink scheduling.
In terms of the classification of \hyperref[sec:taxonomy]{Section~\ref{sec:taxonomy}}, each extension corresponds to a change in one field of the general model, namely in $\gamma$ for objective function variations, in $\beta_3$ for different observation requirements, in $\beta_4$ for image quality requirements, in $\alpha_1$ for active imaging satellites, and in $\beta_1$ for integrated scheduling.
Since not all extensions have been proposed for every formulation category, the coverage of extensions varies across categories. 

\autoref{tab:models} summarizes the most relevant modeling characteristics of the considered formulations. Proposed extensions refer to the additional problem dimensions that have been addressed in the literature for each formulation category. Critical modeling issues refer to the main structural challenges that arise when modeling the EOS scheduling problem. Critical problem features refer to the specific problem characteristics that introduce modeling complexities in each formulation category.

\newcolumntype{P}[1]{>{\raggedright\arraybackslash}p{#1}}

\begin{table}[!htbp]
\caption{General characteristics of the considered optimization models.}
\label{tab:models}
\centering
\scriptsize
\renewcommand{\arraystretch}{1.3}
\begin{tabular}{P{2cm}|P{2.5cm}P{2.5cm}|P{2.5cm}|P{2.5cm}P{2.5cm}}
\hline
\textbf{Time representation}
  & \multicolumn{2}{c|}{Discrete local time}
  & \multicolumn{1}{c|}{Discrete global time}
  & \multicolumn{2}{c}{Continuous time} \\
\hline
{Event representation}
  & OTW selection (\hyperref[sec:time_slots]{Section~\ref{sec:time_slots}})
  & Arc flow (\hyperref[sec:arc_flow_local]{Section~\ref{sec:arc_flow_local}})
  & Time indexed (\hyperref[sec:time_indexed]{Section~\ref{sec:time_indexed}})
  & Disjuntive (\hyperref[sec:time_decision]{Section~\ref{sec:time_decision}})
  & Arc flow (\hyperref[sec:arc_flow_continuous]{Section~\ref{sec:arc_flow_continuous}}) \\
\hline
\textbf{Main decisions}
& \multicolumn{2}{c|}{Allocation, sequencing} 
  & Allocation, sequencing, timing
  & \multicolumn{2}{c}{Allocation, sequencing, timing} \\
\hline
\multirow[t]{6}{2.5cm}{{Key variables}}
  & {\textbf{allocation} (discrete)}
  & {\textbf{allocation} (discrete)}
  & {\textbf{allocation} (discrete)}
  & {\textbf{allocation} (discrete)}
  & {\textbf{allocation} (discrete)} \\
  & $\text{x}_{ijkm}$ defines if satellite $i$ observes strip $j$ in orbit $k$ in OTW $m$
  & $\text{x}_{ijj'kmm'}$ defines if satellite $i$ observes strip $j'$ in OTW $m'$ in orbit $k$ immediately after strip $j$ in OTW $m$
  & $\text{x}_{ijkt}$ defines if satellite $i$ observes strip $j$ in orbit $k$ at time point $t$
  & $\text{x}_{ijk}$ defines if satellite $i$ observes strip $j$ in orbit $k$
  & $\text{x}_{ijj'k}$ defines if satellite $i$ observes strip $j'$ in orbit $k$ immediately after strip $j$ \\
  & {\textbf{sequencing} (discrete)}
  &
  & {\textbf{sequencing} (discrete)}
  & {\textbf{sequencing} (discrete)}
  & {\textbf{timing} (continuous)} \\
  & $\text{y}_{ijj'kmm'}$ defines if satellite $i$ observes strip $j'$ in orbit $k$ in OTW $m'$ immediately after strip $j$ in OTW $m$
  &
  & \textit{general (immediate) precedence:} $\text{y}_{ijj'k}$ defines if satellite $i$ observes strip $j'$ in orbit $k$ (immediately) after strip $j$
  & \textit{general (immediate) precedence:} $\text{y}_{ijj'k}$ defines if satellite $i$ observes strip $j'$ in orbit $k$ (immediately) after strip $j$
  & $\text{t}_{ijk}$ defines the observation start time of strip $j$ by satellite $i$ in orbit $k$ \\
  &
  &
  &
  & {\textbf{timing} (continuous)} \\
  &
  &
  &
  & $\text{t}_{ijk}$ defines the observation start time of strip $j$ by satellite $i$ in orbit $k$
  &  \\
\hline
Proposed extensions
  & Profit based on coverage; Time-dependent profit; Multiple observations
  & Profit based on number of observations; Multiple observations
  & Time-dependent profit; integrated scheduling
  & Active imaging; integrated scheduling; image quality constraints; multiple observations; stereoscopic requests; scheduling without orbit decomposition; time-dependent profit; multi-objective \newline
  & Profit based on observation interval; integrated scheduling \\
Critical modeling issues
  & Binary decision variables scale with the 
  number of OTWs per VTW; binary auxiliary variables and sequencing constraints grow quadratically with OTWs; pre-computation cost of transition times; time interval size must be tuned \newline
  & Binary decision variables grow quadratically with the number of OTWs; pre-computation cost of transition times; time interval size must be tuned
  & Binary decision variables scale with the total number of time intervals; big-${M}$ required for sequencing; pre-computation cost of transition times; time interval size must be tuned
  & Big-${M}$ required for sequencing; nonlinear constraints under time-dependent transitions and image quality
  & Binary decision variables grow quadratically with the number of strips; big-${M}$ required for sequencing; nonlinear constraints under time-dependent transitions and image quality\\
Critical problem features
  & Active imaging; integrated scheduling
  & Active imaging; integrated scheduling
  & Sequence and time-dependent transition times
  & Sequence and time-dependent transition times; time-dependent image quality
  & Sequence and time-dependent transition times; time-dependent image quality \\
\hline
\end{tabular}
\end{table}

\subsection{Local time discretization}
\label{sec:local_time_disc}

The local time discretization divides each VTW into a set of discrete OTWs, each representing a feasible observation opportunity, based on a fixed time interval size and the observation duration of each strip \citep{Gabrel1997,Wang2016,Wang2019b}. That is, each VTW $vw_{ijk}$ is divided into a set of OTWs $\mathbf{OW}_{ijk}$, where each OTW $m \in \mathbf{OW}_{ijk}$, denoted as $ow_{ijkm} = [ows_{ijkm}, owe_{ijkm}]$, is defined by a fixed start time $ows_{ijkm}$ and end time $owe_{ijkm}$, with $owe_{ijkm} - ows_{ijkm} = dur^o_{j}$. 
Since decision variables are defined only over these predefined
opportunities, the observation start time is fixed once an OTW is
selected, and constraint \hyperref[tab:constraints]{I3} is satisfied by
construction. Likewise, any time-dependent parameter can be
pre-computed, namely the transition times between OTWs assigned to the
same satellite, so that $\text{trans}(ows_{ijkm} + dur^o_{j},
ows_{ij'km'}) = \text{trans}^{o2o}_{ijj'kmm'}$. This eliminates the nonlinearities that would otherwise arise from the dependence of transition times on continuous decision variables. The same applies to any other time-dependent parameter incorporated into the model, such as image quality functions and time-dependent profit definitions. 
It is worth noting, however, that while pre-computation eliminates nonlinearities within the optimization model, it introduces a preprocessing cost that scales quadratically with the number of OTWs per orbit, since transition times must be evaluated for every pair of candidate tasks assigned to the same satellite. 
Furthermore, let, for each $i \in \mathbf{I}$, $k \in \mathbf{K}_i$, each pair of tasks $j, j' \in \mathbf{J}$ and each $m \in \mathbf{OW}_{ijk}$, 
$
  \mathbf{OW}^{+}_{ijj'k}(m) = \bigl\{\, m' \in \mathbf{OW}_{ij'k}
  \;:\; ows_{ij'km'} \geq owe_{ijkm}
  + \text{trans}^{o2o}_{ijj'kmm'} \,\bigr\}
$
be the set of OTWs of task $j'$ that can be observed immediately after OTW $ow_{ijkm}$, and let $\mathbf{OW}^{-}_{ijj'k}(m')$ be defined analogously as the set of OTWs of task $j$ that can be observed immediately before OTW $ow_{ij'km'}$. Since the start times of all OTWs are known prior to solving, both sets can also be pre-computed.

\subsubsection{OTW selection}
\label{sec:time_slots}

The most relevant contributions for this representation are \citet{Wang2011,Wang2019b}. In the OTW selection formulation, each OTW is treated as a candidate observation task, and the scheduling problem reduces to selecting a feasible subset of these candidates. Since transition times between all pairs of OTWs are known prior to solving, the feasibility of each pair can be determined in advance and encoded in a binary parameter $\delta_{ijj'kmm'} \in \{0,1\}$, equal to 1 if OTWs $ow_{ijkm}$ and $ow_{ij'km'}$ can both be assigned to satellite $i$ in orbit $k$, that is, if the transition time between them is satisfied in the temporal order induced by their fixed start times, and 0 otherwise. This parameter is defined for every pair of OTWs assigned to the same satellite and orbit, so its number of entries grows quadratically with the number of candidate tasks per orbit, as does the number of sequencing constraints, discussed below.
 The main decision variable is defined as $\text{x}_{ijkm} \in \{0,1\}$, equal to 1 if OTW $ow_{ijkm}$ is selected, and 0 otherwise.

\subsubsubsection{Allocation constraints}
\noindent Constraints~\eqref{eq:slot_selection_I2} ensure that at most one OTW is selected per strip (\hyperref[tab:constraints]{I2}). 
\begin{align}
    & \sum_{i \in \mathbf{I}} \sum_{k \in \mathbf{K}_{ij}} 
    \sum_{m \in \mathbf{OW}_{ijk}} \text{x}_{ijkm} \leq 1, 
    \quad \forall \ \ j \in \mathbf{J} 
    \label{eq:slot_selection_I2}
\end{align}

\subsubsubsection{Sequencing constraints}
\noindent Constraints \eqref{eq:slot_selection_T} make use of the pre-computed parameter $\delta_{ijj'kmm'}$ to ensure that if two OTWs are selected for the same satellite then the transition time between them is feasible (\hyperref[tab:constraints]{I1,I4}).
\begin{align}
    & \text{x}_{ijkm} + \text{x}_{ij'km'} 
    \leq 1 + \delta_{ijj'kmm'}, \quad \forall \ \ i \in \mathbf{I},\; 
    j, j' \in \mathbf{J}\;(j \neq j'),\; k \in \mathbf{K}_i,\; 
    m \in \mathbf{OW}_{ijk},\; m' \in \mathbf{OW}_{ij'k} 
    \label{eq:slot_selection_T}
\end{align}
Note that these constraints do not identify which OTW immediately precedes another: they are pairwise compatibility constraints imposed over all pairs assigned to the same satellite and orbit. Since the temporal order of any two OTWs is fixed by their pre-computed start times, enforcing pairwise feasibility over all pairs implies that the transition time is respected between the OTWs that end up consecutive in the solution, without the need to model succession explicitly.

\subsubsubsection{Resource constraints}
\noindent Since the energy consumed during attitude maneuvers between two consecutive 
observations depends on the selected OTWs, the
auxiliary binary variable $\text{y}_{ijj'kmm'} \in \{0,1\}$ is introduced, equal to 1 
if OTW $ow_{ijkm}$ is immediately followed by OTW $ow_{ij'km'}$, and 0 
otherwise. This variable is defined only for pairs of OTWs that can be executed consecutively, that is, for $m' \in \mathbf{OW}^{+}_{ijj'k}(m)$.
Constraints~\eqref{eq:slot_selection_energy} ensure that the total 
energy consumed by observations and attitude maneuvers in each orbit does 
not exceed the available energy budget (\hyperref[tab:constraints]{R1}), and 
constraints~\eqref{eq:slot_selection_memory} ensure that the total memory 
consumed by observations in each orbit does not exceed the available memory 
capacity~(\hyperref[tab:constraints]{R2}) \citep{Wang2019b}.
\begin{align}
    & \sum_{j \in \mathbf{J}} \sum_{m \in \mathbf{OW}_{ijk}} 
    \bar{E}_i^{o} \cdot dur^o_{j} \cdot \text{x}_{ijkm} 
    +  \sum_{j,j' \in \mathbf{J} (j \neq j')} \sum_{m \in \mathbf{OW}_{ijk}} \sum_{m' \in \mathbf{OW}^{+}_{ijj'k}(m)} 
    \bar{E}_i^{a} \cdot \text{trans}^{o2o}_{ijj'kmm'} \cdot \text{y}_{ijj'kmm'} 
    \leq E^{\max}_i, \label{eq:slot_selection_energy}\\ &\quad \forall \ \ i \in \mathbf{I},\ k \in \mathbf{K}_i \nonumber\\
    & \sum_{j \in \mathbf{J}} \sum_{m \in \mathbf{OW}_{ijk}} \bar{M}_i^{o} 
    \cdot dur^o_{j} \cdot \text{x}_{ijkm} \leq M^{\max}_i, \quad \forall \ \ i \in \mathbf{I},\ k \in \mathbf{K}_i
    \label{eq:slot_selection_memory}
\end{align}
For the consistent definition of the auxiliary variables, for each
satellite $i$ and orbit $k$ a dummy source task $j^s_{ik}$ and a dummy
sink task $j^e_{ik}$ are introduced, each with a single OTW
($m^s_{ik}$ and $m^e_{ik}$, respectively) of zero duration and with
zero transition times to and from any other OTW, so that all arcs
incident to them are feasible by construction. The sets $\mathbf{OW}^{+}_{ijj'k}(m)$ and $\mathbf{OW}^{-}_{ijj'k}(m')$ are extended accordingly to the dummy tasks. The following flow conservation constraints are then required. Constraints~\eqref{eq:slot_source_sink} ensure that the source 
and sink tasks are always selected. Constraints~\eqref{eq:slot_flow_out} 
and~\eqref{eq:slot_flow_in} ensure that each selected OTW has exactly one 
immediate successor and predecessor, respectively. 
\begin{align}
    & \text{x}_{ij^s_{ik}km^s_{ik}} = \text{x}_{ij^e_{ik}km^e_{ik}} = 1,
    \quad \forall\, i \in \mathbf{I},\; k \in \mathbf{K}_i
    \label{eq:slot_source_sink}\\[4pt]
    & \sum_{j' \in \mathbf{J} \cup \{j^e_{ik}\}}
      \sum_{m' \in \mathbf{OW}^{+}_{ijj'k}(m)}
      \text{y}_{ijj'kmm'} = \text{x}_{ijkm},
    \quad \forall\, i \in \mathbf{I},\; k \in \mathbf{K}_i,\;
    j \in \mathbf{J} \cup \{j^s_{ik}\},\; m \in \mathbf{OW}_{ijk}
    \label{eq:slot_flow_out}\\[4pt]
    & \sum_{j \in \mathbf{J} \cup \{j^s_{ik}\}}
      \sum_{m \in \mathbf{OW}^{-}_{ijj'k}(m')}
      \text{y}_{ijj'kmm'} = \text{x}_{ij'km'},
    \quad \forall\, i \in \mathbf{I},\; k \in \mathbf{K}_i,\;
    j' \in \mathbf{J} \cup \{j^e_{ik}\},\; m' \in \mathbf{OW}_{ij'k}
    \label{eq:slot_flow_in}
\end{align}

\subsubsubsection{Extensions}
\noindent \textit{a. Objective function variations}

Here we discuss the different profit definitions considered in the literature for this formulation. \citet{Wang2011} considered spot and polygon targets, and modeled profit based on the percentage of observed area of each target using the piecewise linear function mentioned in \hyperref[sec:object_function]{Section \ref{sec:object_function}}.
\citet{Wang2019b} defined profit as a function of the observation angles, considering that observations closer to the nadir point give more imaging quality and therefore higher profit (\hyperref[tab:constraints]{I5}). The authors defined an adjusted weight for each OTW reflecting this property.
Actually, the local time discretization offers a structural advantage over continuous 
formulations in this definition of profit. In the latter, image
quality depends on the observation angle at the moment of acquisition,
and therefore on the start time, making the objective nonlinear. Here,
the observation angles of each OTW are known in advance, so the
adjusted profit is itself a fixed parameter, and time-dependent profit
reduces to a weighted objective at no additional modeling cost.

\noindent \textit{b. Different observation requirements}

\citet{Wang2019b} considered multiple visit 
requirements, where each strip may be observed more than once, including in the same VTW, being one of the few works to consider this explicitly. In this formulation, the 
adaptation is straightforward, with the allocation
constraint \eqref{eq:slot_selection_I2} being relaxed to allow each strip
to be observed up to $N_j$ times, and the sequencing constraints
\eqref{eq:slot_selection_T} extended to pairs of OTWs of the same
strip, so that transition feasibility is also enforced between
successive observations of a single target. 
Notably, the local time discretization is the only time representation that 
accommodates multiple observations of the same strip within the same 
VTW in a natural way, as the additional observation 
candidates are already encoded as distinct OTWs. Other event 
representations could in principle support this extension, but would 
require the introduction of additional variables and constraints.

Finally, we note that this formulation has not yet been extended 
to address the integrated scheduling problem, nor to account for the 
characteristics of active imaging satellites. We note that
\citet{Wang2011} considers the integrated problem, but the downlink 
tasks are modeled with a different time representation from the 
one described here, and therefore the corresponding formulation is 
not addressed in this subsection.
Extending the OTW selection formulation to the integrated problem 
would require discretizing communication windows into candidate 
tasks similarly to the construction of OTWs, introducing both 
types of candidates into the same scheduling problem. This would 
require redefining the feasibility parameter 
$\delta_{ijj'kmm'}$ for all pairs across both task types, and 
since the sequencing constraints and auxiliary variables already 
grow quadratically with the number of OTWs per orbit, applying 
this structure to the enlarged combined candidate set would 
further amplify model scale. Regarding active imaging satellites, 
a relevant modeling challenge arises from the possibility of 
variable observation durations. In the current formulation, the 
observation duration is a fixed parameter, so all OTWs within a 
given VTW share the same duration. If multiple durations are 
feasible for the same strip and VTW, each combination of start 
time and duration would need to be enumerated as a distinct OTW, 
multiplying the number of candidate tasks per VTW by the number 
of feasible durations and thereby compounding the scalability 
challenges already inherent to this formulation.

\subsubsection{Arc flow}
\label{sec:arc_flow_local}
Another formulation based on local time discretization is the 
arc flow formulation, with \citet{Gabrel1997, Wang2021b,Wang2025b} 
as the main references. This formulation shares the same discretization approach as the OTW
selection formulation, but models the scheduling problem for each
satellite $i$ and orbit $k$ as a path-selection problem on a directed
graph, in which OTWs are represented as nodes and a directed arc
between two nodes is selected if one OTW is the immediate successor of
the other. The main decision variable is 
therefore defined as $\text{x}_{ijj'kmm'} \in \{0,1\}$, equal to 
1 if OTW $ow_{ij'km'}$ is the immediate successor 
of OTW $ow_{ijkm}$ for satellite $i$ in orbit $k$, and 0 
otherwise. With this representation, and because the transition times can be 
computed in advance, decision variables are only defined for pairs of OTWs 
that respect the transition time feasibility, that is, for 
$m' \in \mathbf{OW}^{+}_{ijj'k}(m)$. 
For the consistent definition of the variables, we 
introduce for each satellite $i$ and orbit $k$ a dummy source task 
$j^s_{ik}$ and a dummy sink task $j^e_{ik}$, each with a single 
OTW ($m^s_{ik}$ and $m^e_{ik}$, respectively) of zero duration and with 
zero transition times to and from any other OTW, so that all arcs incident 
to them are feasible by construction. 
Since this is an arc flow formulation and only feasible transitions are 
considered, sequencing constraints \hyperref[tab:constraints]{I1,I4} are 
inherently respected.

\subsubsubsection{Allocation constraints}
\noindent Constraints~\eqref{eq:arcflow_I2} ensure that each strip is 
observed at most once (\hyperref[tab:constraints]{I2}). Constraints~\eqref{eq:arcflow_source} 
and~\eqref{eq:arcflow_flow} are the flow conservation constraints, where constraints \eqref{eq:arcflow_source} ensure that the source and sink tasks are always selected, and constraints~\eqref{eq:arcflow_flow} ensure that each selected OTW has exactly one immediate successor and predecessor.
\begin{align}
    & \sum_{i \in \mathbf{I}} \sum_{k \in \mathbf{K}_{ij}}
    \sum_{m \in \mathbf{OW}_{ijk}}
    \sum_{j' \in \mathbf{J} \cup \{j^e_{ik}\}} 
    \sum_{m' \in \mathbf{OW}^{+}_{ijj'k}(m)}
    \text{x}_{ijj'kmm'} \leq 1,
    \quad \forall \ \ j \in \mathbf{J}
    \label{eq:arcflow_I2}\\[4pt]
    & \sum_{j \in \mathbf{J} \cup \{j^e_{ik}\}} 
    \sum_{m \in \mathbf{OW}^{+}_{ij^s_{ik}jk}(m^s_{ik})}
    \text{x}_{ij^s_{ik}jkm^s_{ik}m} 
    = \sum_{j \in \mathbf{J} \cup \{j^s_{ik}\}} 
    \sum_{m \in \mathbf{OW}^{-}_{ijj^e_{ik}k}(m^e_{ik})}
    \text{x}_{ijj^e_{ik}kmm^e_{ik}} = 1,
    \quad \forall \ \ i \in \mathbf{I},\ k \in \mathbf{K}_i
    \label{eq:arcflow_source} \\[4pt]
    & \sum_{j' \in \mathbf{J} \cup \{j^e_{ik}\}} 
    \sum_{m' \in \mathbf{OW}^{+}_{ijj'k}(m)}
    \text{x}_{ijj'kmm'}
    =
    \sum_{j' \in \mathbf{J} \cup \{j^s_{ik}\}} 
    \sum_{m' \in \mathbf{OW}^{-}_{ij'jk}(m)}
    \text{x}_{ij'jkm'm},
        \label{eq:arcflow_flow}\\
    & 
    \quad \forall \ \ i \in \mathbf{I},\ k \in \mathbf{K}_i,\
    j \in \mathbf{J},\ m \in \mathbf{OW}_{ijk} \nonumber
\end{align}

\subsubsubsection{Resource constraints}
\noindent Constraints~\eqref{eq:arcflow_energy} ensure that the total 
energy consumed by observations and attitude maneuvers in each orbit does 
not exceed the available energy budget~(\hyperref[tab:constraints]{R1}), and 
constraints~\eqref{eq:arcflow_memory} ensure that the total memory 
consumed by observations in each orbit does not exceed the available memory 
capacity~(\hyperref[tab:constraints]{R2}) \citep{Wang2021,Wang2025b}.

It is worth noting that, unlike the OTW selection formulation, 
the arc flow representation does not require auxiliary variables 
to track direct task succession for the purpose of resource 
constraints. Since the arc variable $\text{x}_{ijj'kmm'}$ already 
encodes immediate precedence by definition, the maneuver energy 
term can be expressed directly in terms of the arc variables, 
as shown in constraint~\eqref{eq:arcflow_energy}, without 
introducing additional binary variables.
\begin{align}
& \sum_{j \in \mathbf{J}} \sum_{m \in \mathbf{OW}_{ijk}}
  \sum_{j' \in  \mathbf{J} \cup \{j^e_{ik}\}} 
  \sum_{m' \in \mathbf{OW}^{+}_{ijj'k}(m)}
  \bar{E}_i^{o} \cdot dur^o_{j} \cdot \text{x}_{ijj'kmm'} + \nonumber \\
& \quad + \sum_{j,j' \in \mathbf{J} (j \neq j')} \sum_{m \in \mathbf{OW}_{ijk}}
   \sum_{m' \in \mathbf{OW}^{+}_{ijj'k}(m)}
  \bar{E}_i^{a} \cdot \text{trans}^{o2o}_{ijj'kmm'} \cdot \text{x}_{ijj'kmm'}
\leq E^{\max}_i, \quad \forall\, i \in \mathbf{I},\ k \in \mathbf{K}_i
     \label{eq:arcflow_energy}\\
    & \sum_{j \in \mathbf{J}} \sum_{m \in \mathbf{OW}_{ijk}}
    \sum_{j' \in  \mathbf{J} \cup \{j^e_{ik}\}} 
    \sum_{m' \in \mathbf{OW}^{+}_{ijj'k}(m)}
    \bar{M}_i^{o} \cdot dur^o_{j} \cdot
    \text{x}_{ijj'kmm'}
    \leq M^{\max}_i, \quad \forall\, i \in \mathbf{I},\ k \in \mathbf{K}_i
    \label{eq:arcflow_memory}
\end{align}

\subsubsubsection{Extensions}
\noindent \textit{a. Objective function variations}

The consideration of multiple observations of the same target by \cite{Wang2021b} led to a different definition of profit that considers that the profit of each target increases with the number of observations, up to a maximum 
desired value. To model this, the authors introduce a bounded 
continuous profit variable per target, linked to the observation 
decisions by an accumulation constraint, and replace the standard 
weighted sum in the objective by the sum of these variables. 
They then extend their formulation to handle cloud coverage 
uncertainty, under which the profit of each observation is 
not known exactly at scheduling time, as cloud coverage may reduce 
image quality and therefore its value.
\citep{Wang2025b} also model profit as a function of the number of 
observations per target, but adopt a nonlinear piecewise 
structure to capture diminishing marginal returns.

\noindent \textit{b. Different observation requirements}

\citet{Wang2021b,Wang2025b} 
consider multiple observations of the same target, with the modeling being similar to what was described for
\citep{Wang2019b} in the OTW selection formulation. This is, as already discussed, a simple adpatation of the allocation constraint.

Finally, we note that this formulation has not been extended 
to address the integrated imaging and downlink scheduling 
problem, nor to account for the characteristics of active 
imaging satellites. Extending the arc flow formulation to the 
integrated problem would require redefining the directed graph 
to include nodes of different types, which substantially changes 
the path-selection problem that underlies this representation. 
Furthermore, since the binary decision variables
already grow quadratically with the number of OTWs per orbit, 
enlarging the node set to include downlink candidates would 
further amplify model scale. Regarding active imaging satellites, 
the same challenge identified for the OTW selection formulation 
applies. Variable observation durations would require enumerating 
each combination of start time and duration as a distinct OTW, 
increasing the number of nodes in the graph and thereby 
compounding the scalability challenges already inherent to this 
formulation.

\subsection{Global time discretization}
\label{sec:global_time_disc}

The global time discretization divides the planning horizon 
into $T-1$ equal intervals, defining the set of time points as
$\mathbf{T} = \{1, \ldots, T\}$ \citep{Valicka2019}.
As with local time discretization, this approach allows all 
time-dependent parameters to be pre-computed prior to solving. 
The key difference lies in how candidate positions are generated. 
Rather than dividing each VTW independently, 
the global approach defines a single uniform grid over the entire 
scheduling horizon, so that variables are in principle defined for 
all time points, including those that fall outside the visibility 
windows of any given strip.
Regarding transition times, if we denote by $t, t' \in \mathbf{T}$ 
the observation start time points of strips $j$ and $j'$, 
respectively, then $\text{trans}(t + dur_{j}^o,t') = 
\text{trans}^{o2o}_{ijj'ktt'}$ is defined as a fixed input 
parameter for each pair of time points.

\subsubsection{Time indexed}
\label{sec:time_indexed}

The most relevant contributions for this representation are 
\citet{Lin2005}, \citet{Valicka2019}, \citet{Peng2020}, 
\citet{Du2020}, and \citet{He2022}. In the time indexed 
formulation, each time point is a candidate start time for an 
observation, and the scheduling problem reduces to assigning 
each strip to at most one time point within its feasible set. 
The main decision variable is defined as $\text{x}_{ijkt} 
\in \{0,1\}$, equal to 1 if the observation of strip $j$ by 
satellite $i$ in orbit $k$ starts at time point $t \in \mathbf{T}$, and 0 
otherwise. 
Some authors, such as \citet{Peng2020}, define 
the decision variables only for $t \in \mathbf{T}^o_{ijk}$, where
$\mathbf{T}^o_{ijk} = \{t \in \mathbf{T} : vws_{ijk} \leq t 
\leq vwe_{ijk} - dur_{j}^o\}$, which reduces the number of 
variables and eliminates the need for explicit timing 
constraints. In the following, we adopt the former definition for generality.

\subsubsubsection{Allocation constraints}
\noindent Constraints~\eqref{eq:indexed_I2} ensure that each strip is observed at most 
once (\hyperref[tab:constraints]{I2}).
\begin{align}
    \sum_{i \in \mathbf{I}} \sum_{k \in \mathbf{K}_{ij}} \sum_{t \in \mathbf{T}}
    \text{x}_{ijkt} \leq 1,
    \quad \forall \ \ j \in \mathbf{J}
    \label{eq:indexed_I2}
\end{align}

\subsubsubsection{Timing constraints}
\noindent Constraints~\eqref{eq:indexed_tw_lb}--\eqref{eq:indexed_tw_ub} ensure that each observation is scheduled within its corresponding VTW (\hyperref[tab:constraints]{I3}).
\begin{align}
    &\sum_{t \in \mathbf{T}} t \cdot\text{x}_{ijkt}
    \geq vws_{ijk} \cdot \sum_{t \in \mathbf{T}} \text{x}_{ijkt}, \quad \forall \ \ i \in \mathbf{I},\, j \in \mathbf{J},\, k \in \mathbf{K}_{ij}
    \label{eq:indexed_tw_lb}
    \\[4pt]
    &\sum_{t \in \mathbf{T}} t \cdot \text{x}_{ijkt}
    \leq \bigl(vwe_{ijk} - dur_{j}^o\bigr) \cdot
    \sum_{t \in \mathbf{T}} \text{x}_{ijkt}, \quad \forall \ \ i \in \mathbf{I},\, j \in \mathbf{J},\, k \in \mathbf{K}_{ij}
    \label{eq:indexed_tw_ub}
\end{align}
We note that if we define the decision variable $\text{x}_{ijkt}$ only for $t \in \mathbf{T}^o_{ijk}$, constraints~\eqref{eq:indexed_tw_lb}--\eqref{eq:indexed_tw_ub} 
are not required. 

\subsubsubsection{Sequencing constraints}
\noindent To model the transition constraints, binary auxiliary variables must be 
introduced to encode the sequencing of observations. As discussed in 
\hyperref[sec:classification_models]{Section \ref{sec:classification_models}}, sequencing can be based on either immediate or general precedence. 
Both approaches have been proposed in the literature for the time indexed 
formulation and are described below.

Based on immediate precedence, we define the auxiliary 
binary variables $\text{y}_{ijj'k} \in \{0,1\}$, equal to 1 if strip $j'$ is observed immediately after strip $j$ by satellite 
$i$ in orbit $k$, and 0 otherwise \citep{Peng2020}. Based on general precedence, we define the auxiliary 
binary variables $\text{y}_{ijj'k} \in \{0,1\}$, equal to 1 if strip $j'$ is observed after strip $j$ (not necessarily immediately) by satellite 
$i$ in orbit $k$, and 0 otherwise \citep{Du2020,He2022}. 
Although sequence and time-dependent transition times 
(\hyperref[tab:constraints]{A1c}) do not introduce 
nonlinearity in the time indexed formulation, as transition 
times can always be pre-computed, they 
do introduce a different form of complexity. Unlike 
sequence-dependent transition times (\hyperref[tab:constraints]{A1b}), 
which require one sequencing constraint per ordered pair of 
strips, time-dependent transition times additionally require one constraint 
per pair of time points, increasing the number of sequencing 
constraints and potentially compromising computational efficiency.
Constraints~\eqref{eq:indexed_seq_td} are defined for sequence and time-dependent transition times \citep{Du2020} and constraints~\eqref{eq:indexed_seq_sd} are defined for sequence-dependent transition times \citep{He2022}, to ensure that between two observations of the same satellite the corresponding transition time is respected (\hyperref[tab:constraints]{I1,I4}), where $\mathcal{M}$ is defined as a sufficiently large positive constant. 
\begin{subequations}
    \label{eq:indexed_seq}
\begin{align}
            & t + dur_{j}^o + \text{trans}^{o2o}_{ijj'ktt'}
    \leq t' 
    + \mathcal{M}(3 - \text{x}_{ijkt} - \text{x}_{ij'kt'} - \text{y}_{ijj'k}),
    \label{eq:indexed_seq_td}\\
    & \quad \forall \ \ i \in \mathbf{I},\;
    j, j' \in \mathbf{J} \ (j \neq j'),\;
    k \in \mathbf{K}_i,\;
    t, t' \in \mathbf{T},\ \nonumber \\
    &\sum_{t \in \mathbf{T}} \bigl(t + dur_{j}^o + \text{trans}^{o2o}_{ijj'k}\bigr) \cdot\,
    \text{x}_{ijkt}
    \;\leq\;
    \sum_{t' \in \mathbf{T}} t'\cdot \text{x}_{ij'kt'}
    + \mathcal{M}\bigl(1 - \text{y}_{ijj'k}\bigr),
    \label{eq:indexed_seq_sd}\\
    & \quad \forall \ \ i \in \mathbf{I},\;
    j, j' \in \mathbf{J} \ (j \neq j'),\;
    k \in \mathbf{K}_i \nonumber
\end{align}
\end{subequations}
For the consistent definition of the auxiliary variables based on immediate precedence, a dummy source task $j^s_{ik}$ and dummy sink task $j^e_{ik}$ for 
each satellite and orbit are introduced, with zero duration and zero 
transition times to and from any other task, so that all transitions involving 
them are feasible by construction. 
The following flow conservation constraints are then required. Constraints~\eqref{eq:indexed_flow} ensure that each selected observation has exactly one immediate successor and predecessor, and constraints~\eqref{eq:indexed_source_sink} ensure that the source and sink tasks are always selected.
\begin{align}
    &\sum_{j' \in \mathbf{J} \cup \{j^e_{ik}\} \setminus \{j\}} 
    \text{y}_{ijj'k}
    =
    \sum_{j' \in \mathbf{J} \cup \{j^s_{ik}\} \setminus \{j\}} 
    \text{y}_{ij'jk}
    = \sum_{t \in \mathbf{T}} \text{x}_{ijkt},
    \quad \forall \ \ i \in \mathbf{I},\, j \in \mathbf{J},\, 
    k \in \mathbf{K}_i
    \label{eq:indexed_flow}
    \\[6pt]
    & \sum_{j \in \mathbf{J} \cup \{j^e_{ik}\}} 
    \text{y}_{ij^s_{ik}jk}
    =
    \sum_{j \in \mathbf{J} \cup \{j^s_{ik}\}} 
    \text{y}_{ijj^e_{ik}k}
    = 1,
    \quad \forall \ \ i \in \mathbf{I},\, k \in \mathbf{K}_i
    \label{eq:indexed_source_sink}
\end{align}
For the consistent definition of the auxiliary variables based 
on general precedence,
constraints~\eqref{eq:indexed_anti} are required to ensure 
that exactly one ordering is imposed for each pair of strips.
\begin{align}
    &\text{y}_{ijj'k} + \text{y}_{ij'jk} = 1,
    \quad \forall\ \ i \in \mathbf{I},\ j,j' \in \mathbf{J}\ 
    (j \neq j'),\ k \in \mathbf{K}_i
    \label{eq:indexed_anti}
\end{align}
We note that this form imposes an ordering for every pair of 
strips regardless of whether both are assigned. A tighter 
formulation would condition the constraint on both strips being 
scheduled, forcing both variables to zero whenever at least one 
is unscheduled, which produces a tighter linear programming 
relaxation without affecting the optimal integer solution. 

Therefore, sequencing constraints based on immediate precedence 
are defined by constraints~\eqref{eq:indexed_seq}--\eqref{eq:indexed_source_sink}, 
while sequencing constraints based on general precedence are 
defined by constraints~\eqref{eq:indexed_seq} 
and~\eqref{eq:indexed_anti}.

\subsubsubsection{Resource constraints}
\noindent From the main references, only \citet{Du2020} modeled resource 
constraints within the time indexed formulation. Rather than 
tracking resource consumption as a function of observation 
duration and attitude maneuvers, the authors adopt a simpler 
approach in which each observation of strip $j$ incurs a fixed 
energy cost ${E}_j^{o}$ and occupies a fixed memory volume 
${M}_j^{o}$. This differs from the resource modeling adopted 
in the OTW selection and arc flow formulations, where energy 
consumption depends on both observation duration and attitude 
maneuver costs between consecutive tasks. By treating energy 
and memory as fixed per-observation costs, the formulation 
avoids the need for sequencing variables in the resource 
constraints, at the cost of a less accurate representation 
of the actual onboard resource consumption.
Constraints~\eqref{eq:indexed_energy} ensure that the total 
energy consumed by observations in each orbit does not exceed 
the onboard energy budget (\hyperref[tab:constraints]{R1}), and 
constraints~\eqref{eq:indexed_memory} ensure that the total 
memory occupied does not exceed the onboard storage 
capacity (\hyperref[tab:constraints]{R2}).
\begin{align}
    &\sum_{j \in \mathbf{J}} \sum_{t \in \mathbf{T}}
    {E}_j^{o} \cdot \text{x}_{ijkt}
    \leq E^{\max}_i, \quad \forall \ \ i \in \mathbf{I},\ k \in \mathbf{K}_i
    \label{eq:indexed_energy}
    \\[6pt]
    &\sum_{j \in \mathbf{J}} \sum_{t \in \mathbf{T}}
    {M}_j^{o} \cdot \text{x}_{ijkt}
    \leq M^{\max}_i, \quad \forall \ \ i \in \mathbf{I},\ k \in \mathbf{K}_i
    \label{eq:indexed_memory}
\end{align}

\subsubsubsection{Extensions}
\noindent \textit{a. Objective function variations}

\citet{Peng2020} and \citet{Lin2005} consider time-dependent profits 
in their formulations (\hyperref[tab:constraints]{I5}), where the profit of each observation depends on 
the start time within the VTW. In the time indexed formulation, since the observation angles and therefore 
the profit at each time point are fully determined during preprocessing, 
the time-dependent profit simply becomes a time indexed coefficient in 
the objective function. This is a structural advantage of this 
representation over continuous formulations, where the profit would 
depend on a continuous start time variable, introducing nonlinearity 
into the objective.

\noindent \textit{b. Integrated scheduling problem}

\citet{He2022} extended the time indexed formulation to the integrated
problem. For consistency with the notation
adopted throughout, we present the formulation with separate
allocation variables for observation and downlink tasks.
\citet{He2022} instead employ a single allocation variable whose
interpretation depends on whether the associated time window is a VTW
or a CTW, with the satellite and ground station implicit in the window
index. The two formulations are equivalent. In our notation, to model downlink tasks, a binary decision variable
$\text{z}_{ijklt} \in \{0,1\}$ is introduced, equal to 1 if satellite
$i$ starts downloading the data of strip $j$ to ground station $l$ in orbit $k$
at time point $t$, and 0 otherwise. 
This implies new timing constraints \eqref{eq:indexed_int_tw_lb} and \eqref{eq:indexed_int_tw_ub} analogous to 
\eqref{eq:indexed_tw_lb} and \eqref{eq:indexed_tw_ub}, respectively, to ensure that each downlink must start within the 
corresponding CTW (\hyperref[tab:constraints]{D3}).
\begin{align}
    &\sum_{t \in \mathbf{T}} t \cdot\text{z}_{ijklt}
    \geq cws_{ikl} \cdot \sum_{t \in \mathbf{T}} \text{z}_{ijklt}, \quad \forall \ \ i \in \mathbf{I},\, j \in \mathbf{J},\, l \in \mathbf{L},\, k \in \mathbf{K}_{il}
    \label{eq:indexed_int_tw_lb}
    \\[4pt]
    &\sum_{t \in \mathbf{T}} t \cdot \text{z}_{ijklt}
    \leq \bigl(cwe_{ikl} - dur_{j}^d\bigr) \cdot
    \sum_{t \in \mathbf{T}} \text{z}_{ijklt}, \quad \forall \ \ i \in \mathbf{I},\, j \in \mathbf{J},\, l \in \mathbf{L},\, k \in \mathbf{K}_{il}
    \label{eq:indexed_int_tw_ub}
\end{align}
The observation and downlink decisions must be coupled so that every observed 
strip is downlinked (\hyperref[tab:constraints]{ID1}), which 
constraints~\eqref{eq:indexed_obs_dwn} enforce.
\begin{align}
    & \sum_{l \in \mathbf{L}} \sum_{k \in \mathbf{K}_{il}} \sum_{t \in \mathbf{T}}
    \text{z}_{ijklt}
    = \sum_{k \in \mathbf{K}_{ij}} \sum_{t \in \mathbf{T}} \text{x}_{ijkt},
    \quad \forall \ \ i \in \mathbf{I},\, j \in \mathbf{J}
    \label{eq:indexed_obs_dwn}
\end{align}
Since the observation allocation constraints~\eqref{eq:indexed_I2} already bound 
the right-hand side by one, these constraints also ensure that each strip is 
downlinked at most once. 
In addition to this, we must ensure that a strip can only be downloaded after it 
has been observed (\hyperref[tab:constraints]{ID2}). 
Constraints~\eqref{eq:indexed_causality} ensure that the downlink starts only 
after the observation is complete.
\begin{align}
    & \sum_{l \in \mathbf{L}} \sum_{k' \in \mathbf{K}_{il}} \sum_{t \in \mathbf{T}}
    t \cdot \text{z}_{ijk'lt}
    \;\geq\;
    \sum_{k \in \mathbf{K}_{ij}} \sum_{t \in \mathbf{T}}
    \bigl(t + {dur}^o_{j} \bigr) \cdot \text{x}_{ijkt}, \quad  \forall \ \ i \in \mathbf{I},\, j \in \mathbf{J}
    \label{eq:indexed_causality}
\end{align}
Regarding sequencing constraints, \citet{He2022} assumed
sequence-dependent transition times (\hyperref[tab:constraints]{A1b})
between two observations of the same satellite
($\text{trans}^{o2o}_{ijj'}$) (\hyperref[tab:constraints]{I1,I4}), two
downlinks of the same satellite ($\text{trans}^{d2d}_{ijj'}$)
(\hyperref[tab:constraints]{D1,D4}), and an observation of strip $j$
followed by a downlink of strip $j'$ on the same satellite
($\text{trans}^{o2d}_{ijj'}$) (\hyperref[tab:constraints]{ID3}). The
authors additionally consider the rotation time of the receiving
antenna between two consecutive downlinks at the same ground station,
which they also model as sequence-dependent; in line with the
convention adopted in Section~\ref{sec:problem_desc}, we treat it here
as a constant $\text{trans}^{d2d,gs}$
(\hyperref[tab:constraints]{D2,D5}).
To model these, auxiliary binary variables based on general precedence
are introduced for each transition type: $\text{y}^{o2o}_{ijj'} \in
\{0,1\}$, equal to 1 if strip $j$ is observed before strip $j'$ by
satellite $i$, and 0 otherwise, and analogously for
$\text{y}^{d2d}_{ijj'}$, $\text{y}^{o2d}_{ijj'}$, and
$\text{y}^{d2d,gs}_{jj'l}$, the latter ordering the downlinks of
strips $j$ and $j'$ at ground station $l$. Note that $j'$ is not
required to be the immediate successor of $j$.
For each, constraints analogous to those defined in~\eqref{eq:indexed_seq_sd} and \eqref{eq:indexed_anti} are required to ensure that the corresponding transition times are respected and that the ordering is consistent, respectively.
We note that \citet{He2022} did not consider resource constraints, so we do 
not present the corresponding extension.
The proposed integrated formulation thus consists of allocation 
constraints~\eqref{eq:indexed_I2} and~\eqref{eq:indexed_obs_dwn}, timing constraints \eqref{eq:indexed_tw_lb}--\eqref{eq:indexed_tw_ub}, \eqref{eq:indexed_int_tw_lb}--\eqref{eq:indexed_int_tw_ub} and \eqref{eq:indexed_causality}, and 
sequencing constraints for all transition types 
(analogous to constraints~\eqref{eq:indexed_seq_sd} and \eqref{eq:indexed_anti} for all transition types).

Finally, we note that this formulation has not been extended 
to account for the characteristics of active imaging satellites. 

\subsection{Continuous time}
\label{sec:continuous_time}
In continuous-time formulations, the start time of each observation is 
treated as a continuous decision variable, bounded by the continuous 
VTW $vw_{ijk} = [vws_{ijk}, vwe_{ijk}]$. This has a fundamental 
structural consequence. Since transition times between observations 
depend on the exact start times of the tasks and these are decision 
variables, transition times cannot be pre-computed prior to solving, meaning they are functions of the decision variables. This directly introduces nonlinear constraints. The only tractable approach to avoid nonlinearity is to 
ignore time-dependent transition times.

\subsubsection{Disjuntive}
\label{sec:time_decision}
A wide variety of disjunitve continuous-time formulations have been proposed in the literature in the last years. Some of the most relevant contributions in this category are \citep{Lemaitre2002,Liu2017,Chen2019,Peng2019,Chang2021,Zhang2022,Peng2025,Kim2026}.
In the disjuntive formulation, the scheduling problem is decomposed 
into two main decisions, namely which VTW is assigned to each strip 
and at what instant within that VTW the observation starts.
The key decision variables are the binary variable $\text{x}_{ijk} 
\in \{0,1\}$, equal to 1 if strip $j$ is observed by satellite $i$ 
in orbit $k$, and 0 otherwise, and the continuous variable 
$\text{t}_{ijk} \geq 0$ representing the start time of the observation 
of strip $j$ by satellite $i$ in orbit $k$.

\subsubsubsection{Allocation constraints}
\noindent Constraints~\eqref{eq:cont_alloc} ensure that each strip is 
observed at most once (\hyperref[tab:constraints]{I2}).
\begin{align}
    & \sum_{i \in \mathbf{I}} \sum_{k \in \mathbf{K}_{ij}} 
    \text{x}_{ijk} \leq 1, 
    \quad \forall \ \ j \in \mathbf{J} 
    \label{eq:cont_alloc}
\end{align}
\subsubsubsection{Timing constraints}

\noindent Constraints~\eqref{eq:cont_lb} and~\eqref{eq:cont_ub} ensure that each observation is scheduled within its corresponding VTW (\hyperref[tab:constraints]{I3}), where $\mathcal{M}$ is defined as a sufficiently large positive constant.
\begin{align}
    & \text{t}_{ijk} \geq vws_{ijk} \cdot \text{x}_{ijk} , \quad \forall \ \ i \in \mathbf{I},\; j \in \mathbf{J},\; 
    k \in \mathbf{K}_{ij} 
    \label{eq:cont_lb} \\
    & \text{t}_{ijk} + dur^o_{j} \leq vwe_{ijk} + 
    \mathcal{M}(1 - \text{x}_{ijk}), \quad \forall \ \ i \in \mathbf{I},\; j \in \mathbf{J},\; 
    k \in \mathbf{K}_{ij}
    \label{eq:cont_ub}
\end{align}
We note that constraint~\eqref{eq:cont_ub} can equivalently be stated
without the big-$\mathcal{M}$ parameter. For instance, \citet{Kim2026}
define constraints~\eqref{eq:cont_ub_nobigM}, a variable upper bound
in which the bound on $\text{t}_{ijk}$ scales directly with
$\text{x}_{ijk}$. For $\text{x}_{ijk} \in \{0,1\}$ both formulations
admit the same set of feasible schedules, the only difference being
that \eqref{eq:cont_ub_nobigM} additionally forces
$\text{t}_{ijk} = 0$ whenever strip $j$ is not selected, thereby
fixing a variable that is otherwise meaningless. Its linear relaxation
is, however, strictly tighter, since the bound is relaxed
proportionally to $\text{x}_{ijk}$ rather than by $\mathcal{M}$.
\begin{equation}
    \text{t}_{ijk} \leq \left({vwe}_{ijk} - {dur}^{{o}}_{j}\right) 
    \cdot \text{x}_{ijk}, 
    \qquad \forall\ \ i \in \mathbf{I},\; j \in \mathbf{J},\; k \in \mathbf{K}_{ij}
    \label{eq:cont_ub_nobigM}
\end{equation}

\subsubsubsection{Sequencing constraints}
\noindent Again, to model the transition constraints, binary auxiliary variables must be 
introduced to encode the sequencing of observations, which can be based on either immediate or general precedence. 
Both approaches have been proposed in the literature for this
formulation and are described below.

Based on immediate precedence, we define the auxiliary 
binary variables $\text{y}_{ijj'k} \in \{0,1\}$, equal to 1 if strip $j'$ is observed immediately after strip $j$ by satellite 
$i$ in orbit $k$, and 0 otherwise \citep{Liu2017,He2018,Yang2018,Chang2021,Wang2024,Peng2025}. Based on general precedence, we define the auxiliary 
binary variables $\text{y}_{ijj'k} \in \{0,1\}$, equal to 1 if strip $j'$ is observed after strip $j$ (not necessarily immediately) by satellite 
$i$ in orbit $k$, and 0 otherwise \citep{Xiao2019,Kim2026}. In either case, constraints~\eqref{eq:cont_nooverlap} are defined to ensure that between two observations of the same satellite the corresponding transition time is respected (\hyperref[tab:constraints]{I1,I4}).
We note that, since the transition times 
between observations depend on the exact start times of the tasks and 
cannot be pre-computed prior to solving in this formulation type, when considering sequence and time-dependent transition times (\hyperref[tab:constraints]{A1c}), constraints~\eqref{eq:cont_nooverlap} are nonlinear constraints. 
Because of this, most works assume transition times as sequence-dependent only (\hyperref[tab:constraints]{A1b})
\citep{Lemaitre2002,Chen2019,Cho2018,Zhang2022}, in which case the constraints become linear. The following 
constraints are stated for the general case of time-dependent 
transition times.
\begin{align}
    &\text{t}_{ij'k} \geq \text{t}_{ijk} + dur^o_{j} +  \text{trans}(\text{t}_{ijk} + dur^o_{j},\text{t}_{ij'k})
    - \mathcal{M}(1 -\text{y}_{ijj'k} ),\label{eq:cont_nooverlap}  \\ & \quad \forall \ \ i \in \mathbf{I},\; 
    j, j' \in \mathbf{J}\; (j \neq j'),\; k \in \mathbf{K}_i 
    \nonumber
\end{align}
For the consistent definition of the auxiliary variables based on immediate precedence, a dummy source task $j^s_{ik}$ and dummy sink task $j^e_{ik}$ for each 
satellite $i \in \mathbf{I}$ and orbit $k \in \mathbf{K}_i$ are introduced, with zero duration and zero 
transition times to and from any other task, so that all transitions involving 
them are feasible by construction. .
The following flow conservation constraints are then required to
ensure that each selected observation has exactly one immediate successor and predecessor \eqref{eq:cont_flow}, to ensure that the source and sink tasks are always selected \eqref{eq:cont_source_sink} \citep{Wang2024a,Peng2025} .
\begin{align}
    & \sum_{j' \in \mathbf{J} \cup \{j^e_{ik}\} \setminus \{j\}} \text{y}_{ijj'k}
    = \sum_{j' \in \mathbf{J} \cup \{j^s_{ik}\} \setminus \{j\}} \text{y}_{ij'jk}
    = \text{x}_{ijk}, \quad \forall \ \ i \in \mathbf{I},\, j \in \mathbf{J},\, k \in \mathbf{K}_i
    \label{eq:cont_flow} \\[4pt]
    & \sum_{j' \in \mathbf{J} \cup \{j^e_{ik}\}} \text{y}_{i j^s_{ik} j' k}
    = \sum_{j' \in \mathbf{J} \cup \{j^s_{ik}\}} \text{y}_{i j' j^e_{ik} k}
    = 1, \quad \forall \ \ i \in \mathbf{I},\, k \in \mathbf{K}_i
    \label{eq:cont_source_sink}
\end{align}
For the consistent definition of the auxiliary variables based on general precedence, constraints~\eqref{eq:cont_anti} are required to ensure that exactly one ordering is imposed for each pair of strips, when both strips are assigned \citep{Kim2026}. As these constraints involve the product of two binary 
variables, they are nonlinear and must be linearized prior to solving. 
\begin{align}
    &\text{y}_{ijj'k} + \text{y}_{ij'jk} = 
    \text{x}_{ijk}
    \cdot
     \text{x}_{ij'k}, \quad \forall\ \ i \in \mathbf{I},\ j,j' \in \mathbf{J}\ (j \neq j'),\
    k \in \mathbf{K}_i
    \label{eq:cont_anti}
\end{align}

Therefore, sequencing constraints based on immediate precedence are 
defined by constraints~\eqref{eq:cont_nooverlap}--\eqref{eq:cont_source_sink}, 
while sequencing constraints based on general precedence are defined 
by constraints~\eqref{eq:cont_nooverlap} and \eqref{eq:cont_anti}.
We note that general precedence, as discussed in the resource constraints below, does not allow 
maneuver energy to be modeled directly without additional variables.

\subsubsubsection{Resource constraints}
\noindent Several authors have incorporated resource constraints into the 
continuous-time formulation, with different modeling approaches.
Regarding energy consumption, the most common approach considers energy as a function of both 
observation duration and attitude maneuvers between consecutive tasks. 
As already noted, maneuver energy depends on the timing 
of successive observations, and its inclusion therefore requires an 
auxiliary sequencing variable that identifies direct task succession. 
When sequencing variables are defined under immediate precedence, 
this term is well-defined and constraints~\eqref{eq:cont_energya} can 
be naturally stated to ensure that the total energy consumed per orbit 
does not exceed the available budget (\hyperref[tab:constraints]{R1}) \citep{Gu2022,Wu2022,Zhou2023,Wang2024a}.
When sequencing variables are defined under general precedence, 
incorporating maneuver energy in this manner is no longer 
straightforward. Since $\text{y}_{ijj'k} = 1$ for all ordered pairs 
rather than only adjacent ones, the maneuver term would overcount the 
actual energy consumed. Modeling maneuver energy under general precedence would in practice
require replacing the general precedence variables by immediate
precedence ones, which may explain why this approach has not yet been
adopted in such formulations. Accordingly, authors working 
with general-precedence sequencing variables either simplify the energy 
model to account solely for observation duration, yielding a constraint 
analogous to~\eqref{eq:cont_energya} but without the maneuver term 
\citep{Zhang2022}, or ignore energy consumption constraints \citep{Xiao2019,Kim2026}.
\begin{align}
    &\sum_{j \in \mathbf{J}} \bar{E}_i^{o} \cdot dur_{j}^o 
    \cdot \text{x}_{ijk} + \sum_{j, j' \in \mathbf{J} (j \neq j')} 
    \bar{E}_i^{a} \cdot \text{trans}(\text{t}_{ijk} + dur^o_{j},\text{t}_{ij'k}) \cdot \text{y}_{ijj'k}
    \leq E^{\max}_i,\quad \forall \ \ i \in \mathbf{I}, k \in \mathbf{K}_i
    \label{eq:cont_energya}
\end{align}
We note that under the sequence-dependent simplification discussed 
in the sequencing constraints, $\text{trans}(\text{t}_{ijk} + dur^o_{j},\text{t}_{ij'k})$ 
reduces to a fixed parameter $\text{trans}^{o2o}_{ijj'k}$, and 
constraint~\eqref{eq:cont_energya} remains linear. Without this 
simplification the maneuver term is nonlinear.

Memory consumption is typically modeled as a function of observation duration only, defining constraints~\eqref{eq:cont_memorya} to ensure that the total memory consumed per orbit does not exceed the available storage capacity (\hyperref[tab:constraints]{R2}) \citep{Wang2024a}.
\begin{align}
    &\sum_{j \in \mathbf{J}} \bar{M}_i^{o} \cdot dur_{j}^o 
    \cdot \text{x}_{ijk}
    \leq M^{\max}_i, \quad \forall \ \ i \in \mathbf{I}, k \in \mathbf{K}_i
    \label{eq:cont_memorya}
\end{align}

\subsubsubsection{Extensions}

\noindent \textit{a. Objective function variations}

Being the most widely adopted formulation in the literature, the disjuntive formulation also exhibits the greatest 
diversity of objective functions. Most works consider 
maximization of a weighted linear profit summed over scheduled 
targets. \citet{Wolfe2000,Peng2019,Peng2025} define time-dependent 
profits based on a defined instant image quality function (\hyperref[tab:constraints]{I5}). 
\citet{Zhou2023} additionally accounts for the effect of 
variable observation duration, decomposing the profit into 
two components: one based on the instant image quality at 
the observation start time and the other based on the duration of the observation, defining a logarithmic function of the duration that reflects the diminishing marginal 
returns associated with longer imaging times. 
\citet{Lu2021,Chang2025a,Chang2025b,Wang2025d} take a 
different approach, defining the profit as the fraction of 
total instantaneous quality available in the VTW that is 
accumulated within the chosen observation window. Unlike the 
preceding formulations, where quality is assessed at the 
observation start time and duration enters the profit as a 
separate component, this definition integrates instantaneous 
quality over the full observation window, making the profit 
a unified function of both start time and duration.
 Under the 
continuous-time formulation, such profit definitions are 
inherently time-dependent and nonlinear, which is the 
primary reason why the optimization methods adopted in 
these works rely on non-exact approaches or simplifying 
assumptions. Some works instead consider minimization of the makespan
\citep{Xiao2019,Qu2023,Qu2025}, that is, the completion time of the
last scheduled observation, which is linearized by an auxiliary
variable bounded below by the completion time of every scheduled
observation, or minimization of energy consumption 
\citep{Wang2014,Chang2021,Lu2021,Chang2025a,Chang2025b,Wang2025d}, 
which is nonlinear when transition times are time-dependent and 
linear otherwise. Multi-objective 
formulations are also common, the most frequent combination 
being time-dependent profit and energy consumption 
\citep{Lu2021,Chang2025a,Chang2025b,Wang2025d}, with some 
works additionally maximizing the number of scheduled tasks.

\noindent \textit{b. Integrated scheduling problem}

The extension of continuous-time formulations to the 
integrated scheduling problem has received some attention in the literature 
\citep{Xiao2019,Chang2021,Zhang2022,Kim2026}, 
but the existing contributions differ substantially in 
their structural assumptions, making a single 
canonical formulation difficult to 
establish. We therefore describe the modeling elements 
shared by most contributions and highlight the key 
differentiating assumptions.
To model downlink tasks, a binary variable $\text{z}_{ijkl} \in \{0,1\}$ is 
introduced, equal to 1 if satellite $i$ downloads strip $j$ to ground 
station $l$ in orbit $k$, as well as a continuous decision variable 
$\text{t}^d_{ijkl} \geq 0$ representing the corresponding downlink start time.
This implies new timing constraints \eqref{eq:cont_int_lb} and \eqref{eq:cont_int_ub} analogous to 
\eqref{eq:cont_lb} and \eqref{eq:cont_ub}, respectively, to ensure that each downlink 
must start within the corresponding CTW (\hyperref[tab:constraints]{D3}).
\begin{align}
    & \text{t}^d_{ijkl} \geq cws_{ikl} \cdot \text{z}_{ijkl}, \quad \forall \ i \in \mathbf{I},\; j \in \mathbf{J},\; 
    k \in \mathbf{K}_{il},\; l \in \mathbf{L} 
    \label{eq:cont_int_lb} \\
    & \text{t}^d_{ijkl} + dur^d_{j} \leq cwe_{ikl} + 
    \mathcal{M}(1 - \text{z}_{ijkl}), \quad \forall \ i \in \mathbf{I},\; j \in \mathbf{J},\; 
    k \in \mathbf{K}_{il},\; l \in \mathbf{L} 
    \label{eq:cont_int_ub} 
\end{align}
In addition to this, we must ensure that a strip can only be downloaded after it has been observed (\hyperref[tab:constraints]{ID2}). Constraints~\eqref{eq:cont_int_causal} enforce the temporal 
ordering between observation and downlink.
\begin{align}
    & \text{t}^d_{ijk'l} \geq \text{t}_{ijk} + dur^o_{j} 
    - \mathcal{M}\bigl(2 - \text{x}_{ijk} - 
    \text{z}_{ijk'l}\bigr)
    \label{eq:cont_int_causal} \\
    & \quad \forall \ \ i \in \mathbf{I},\; j \in \mathbf{J},\;
    k \in \mathbf{K}_{ij},\; k' \in \mathbf{K}_{il},\; 
    l \in \mathbf{L} \nonumber
\end{align}
With respect to the downlink allocation constraints, constraints~\eqref{eq:cont_int_alloc_eq} ensure that every observed 
strip is also downlinked (\hyperref[tab:constraints]{ID1}) 
\citep{Xiao2019,Chang2021,Zhang2022}. Since the observation allocation 
constraints~\eqref{eq:cont_alloc} already bound the right-hand side by one, they 
also ensure that each strip is downlinked at most once.
\begin{align}
    & \sum_{l \in \mathbf{L}} \sum_{k \in \mathbf{K}_{il}} 
    \text{z}_{ijkl} 
    = \sum_{k \in \mathbf{K}_{ij}} \text{x}_{ijk},
    \quad \forall \ i \in \mathbf{I},\; j \in \mathbf{J}
    \label{eq:cont_int_alloc_eq}
\end{align}
Regarding sequencing constraints, there are significant modeling variations 
in the literature. \citet{Xiao2019} do not consider transition times between tasks. The
problem is modeled as a two-stage flow shop, in which observations and
downlinks are treated as separate stages, so that the observation of
one strip and the downlink of another may be executed concurrently by
the same satellite. Non-overlap is therefore enforced only within each
stage, between consecutive observations
(\hyperref[tab:constraints]{I1}) and between consecutive downlinks
(\hyperref[tab:constraints]{D1}) of the same satellite, through the
binary auxiliary variables $\text{y}^{o2o}_{ijj'k} \in \{0,1\}$, equal to 1 
if satellite $i$ observes strip $j$ before strip $j'$ in orbit $k$ (not 
necessarily immediately before), and 0 otherwise; and analogously 
$\text{y}^{d2d}_{ijj'k} \in \{0,1\}$, equal to 1 if the downlink of strip $j$ 
precedes that of strip $j'$ on satellite $i$ in orbit $k$, irrespective of the 
ground stations involved, defining constraints 
\eqref{eq:cont_int_nooverlap_o2o_xiao} and 
\eqref{eq:cont_int_nooverlap_d2d_xiao}, respectively.
\begin{align}
    &\text{t}_{ij'k} \geq \text{t}_{ijk} + dur^{o}_{j}
      - \mathcal{M}\bigl(1 - \text{y}^{o2o}_{ijj'k}\bigr),
      \quad \forall\; i \in \mathbf{I},\; j,j' \in \mathbf{J}
      \;(j \neq j'),\; k \in \mathbf{K}_{i}
    \label{eq:cont_int_nooverlap_o2o_xiao} \\
    &\text{t}^d_{ij'kl'} \geq \text{t}^d_{ijkl} + dur^{d}_{j}
      - \mathcal{M}\bigl(3 - \text{y}^{d2d}_{ijj'k}
      - \text{z}_{ijkl} - \text{z}_{ij'kl'}\bigr),
    \label{eq:cont_int_nooverlap_d2d_xiao} \\
    &\quad \forall\; i \in \mathbf{I},\; j,j' \in \mathbf{J}\;(j \neq j'),\;
      l, l' \in \mathbf{L},\; k \in \mathbf{K}_{i}
      \nonumber
\end{align}
For the consistent definition of the auxiliary variables, constraints analogous 
to~\eqref{eq:cont_anti} are required.
Differently, \citet{Chang2021} consider two types of 
transitions: time-dependent transitions between consecutive 
observations of the same satellite 
$\text{trans}(\text{t}_{ijk}+ dur^o_{j},\text{t}_{ij'k})$ (\hyperref[tab:constraints]{I4} under \hyperref[tab:constraints]{A1c}), 
and sequence-independent transitions between consecutive 
downlinks from different satellites to the same ground 
station $\text{trans}^{d2d,gs}$ (\hyperref[tab:constraints]{D2,D5} under \hyperref[tab:constraints]{A1a}). To model these, we use 
the auxiliary binary variable $\text{y}^{o2o}_{ijj'k} \in 
\{0,1\}$ based on immediate precedence, already defined in 
the imaging problem (here with an additional superscript), 
and introduce a binary auxiliary variable 
$\text{y}^{d2d,gs}_{ii'jj'kl} \in \{0,1\}$, equal to 1 if 
the downlink of strip $j$ by satellite $i$ to ground station 
$l$ in orbit $k$ is immediately followed by the downlink of 
strip $j'$ by satellite $i'$ to the same ground station $l$, 
and 0 otherwise. For each transition type, non-overlap and 
flow conservation constraints analogous to 
\eqref{eq:cont_nooverlap}--\eqref{eq:cont_source_sink} can 
be defined.
\citet{Zhang2022} consider three types of 
sequence-independent transitions (\hyperref[tab:constraints]{A1a}): between consecutive 
observations of the same satellite ($\text{trans}^{o2o}$) (\hyperref[tab:constraints]{I4}), 
between consecutive downlinks of the same satellite to 
different ground stations ($\text{trans}^{d2d}$) (\hyperref[tab:constraints]{D1,D4}), and 
between consecutive downlinks from different satellites to 
the same ground station ($\text{trans}^{d2d,gs}$) (\hyperref[tab:constraints]{D2,D5}). To model 
the first two, a binary auxiliary variable 
$\text{y}^{d2d}_{ijj'k} \in \{0,1\}$ is introduced, equal to 1 if 
the downlink of strip $j$ by satellite $i$ in orbit $k$ 
immediately precedes the downlink of strip $j'$ by the same 
satellite to a different ground station, and 0 otherwise. 
Together with $\text{y}^{o2o}_{ijj'k}$ and 
$\text{y}^{d2d,gs}_{ii'jj'kl}$ already introduced, 
non-overlap constraints analogous to 
\eqref{eq:cont_nooverlap} apply for each transition type, as 
well as constraints analogous 
to~\eqref{eq:cont_anti}.
\citet{Kim2026} additionally consider a transition between an
observation and a downlink of the same satellite
($\text{trans}^{o2d}$) (\hyperref[tab:constraints]{ID3}). To model
this, we introduce the binary auxiliary variables
$\text{y}^{o2d}_{ijj'kl}$ and $\text{y}^{d2o}_{ijj'kl} \in \{0,1\}$,
equal to 1 if the observation of strip $j$ by satellite $i$ in orbit
$k$ precedes the downlink of strip $j'$ by the same satellite to
ground station $l$, and conversely, and 0 otherwise. Unlike the
previous contributions, \citet{Kim2026} adopt general precedence and
impose the ordering of each pair only when both tasks are scheduled,
which corresponds to the tighter form discussed in
\hyperref[sec:time_indexed]{Section \ref{sec:time_indexed}}. Together
with $\text{y}^{o2o}_{ijj'k}$, $\text{y}^{d2d}_{ijj'k}$ and
$\text{y}^{d2d,gs}_{ii'jj'kl}$ already introduced, non-overlap
constraints analogous to \eqref{eq:cont_nooverlap} apply for each
transition type, as well as constraints analogous
to~\eqref{eq:cont_anti}.

Resource constraints have also been extended to account for 
energy and memory consumption associated with downlink tasks. 
For energy, the most common approach is to add a downlink 
energy term to the existing observation energy budget, 
modeled analogously by replacing observation duration and 
energy rate with their downlink counterparts. The resulting 
constraint is structurally identical to those defined for the 
imaging problem, with the addition of this new component 
\citep{Zhang2022}. For memory, onboard data storage increases with observations and is
partially recovered through downlinks. Since it can rise and fall
along the orbit, the capacity limit must hold at every task of the
sequence and not only in aggregate. This is enforced by tracking the
data volume stored onboard immediately before each scheduled task,
using the auxiliary precedence variables to identify the observations
and downlinks that precede it \citep{Chang2021,Zhang2022,Kim2026}.
Constraints~\eqref{eq:mem_def_o} and~\eqref{eq:mem_def_d} define the
continuous auxiliary variables $\text{r}^{M,o}_{ijk} \geq 0$ and
$\text{r}^{M,d}_{ijkl} \geq 0$, representing the data stored onboard
satellite $i$ immediately before observing strip $j$ in orbit $k$, and
immediately before its downlink to ground station $l$, respectively.
Here the precedence variables must be defined under general
precedence, as the sums range over all preceding tasks, and
$\text{y}^{d2o}_{ijj'kl}$ is defined analogously to
$\text{y}^{o2d}_{ijj'kl}$ in the reverse order.
\begin{align}
    &\text{r}^{M,o}_{ijk} =
    \sum_{j' \in \mathbf{J} \setminus \{j\}}
    \bar{M}^{o}_i \cdot dur^{o}_{j'} \cdot \text{y}^{o2o}_{ij'jk}
    - \sum_{j' \in \mathbf{J}} \sum_{l \in \mathbf{L}}
    \bar{M}^{d}_i \cdot dur^{d}_{j'} \cdot \text{y}^{d2o}_{ij'jkl}, \quad \forall\ i \in \mathbf{I},\; j \in \mathbf{J},\;
    k \in \mathbf{K}_i,\; l \in \mathbf{L}
    \label{eq:mem_def_o} \\
    &\text{r}^{M,d}_{ijkl} =
    \sum_{j' \in \mathbf{J}}
    \bar{M}^{o}_i \cdot dur^{o}_{j'} \cdot \text{y}^{o2d}_{ij'jkl}
    - \sum_{j' \in \mathbf{J} \setminus \{j\}}
    \bar{M}^{d}_i \cdot dur^{d}_{j'} \cdot \text{y}^{d2d}_{ij'jk}, \quad \forall\ i \in \mathbf{I},\; j \in \mathbf{J},\;
    k \in \mathbf{K}_i,\; l \in \mathbf{L}
    \label{eq:mem_def_d} 
\end{align}
Constraints~\eqref{eq:mem_bound_o} then ensure that enough free memory
is available to store the observation about to be performed, while
constraints~\eqref{eq:mem_bound_d} ensure that the data volume
transmitted in each downlink does not exceed the data accumulated
onboard prior to that task (\hyperref[tab:constraints]{R2}).
\begin{align}
    &\text{r}^{M,o}_{ijk} + \bar{M}^{o}_i \cdot dur^{o}_{j}
    \cdot \text{x}_{ijk} \leq M^{\max}_i, \quad \forall\ i \in \mathbf{I},\; j \in \mathbf{J},\;
    k \in \mathbf{K}_i,\; l \in \mathbf{L} 
    \label{eq:mem_bound_o} \\
    &\bar{M}^{d}_i \cdot dur^{d}_{j} \cdot \text{z}_{ijkl}
    \leq \text{r}^{M,d}_{ijkl}, \quad \forall\ i \in \mathbf{I},\; j \in \mathbf{J},\;
    k \in \mathbf{K}_i,\; l \in \mathbf{L} 
    \label{eq:mem_bound_d} 
\end{align}

The integrated formulation thus consists of 
allocation constraints~\eqref{eq:cont_alloc} 
and~\eqref{eq:cont_int_alloc_eq}, 
timing constraints~\eqref{eq:cont_lb}--\eqref{eq:cont_ub} 
and~\eqref{eq:cont_int_lb}--\eqref{eq:cont_int_causal}, 
sequencing constraints for all transition types considered 
with respective consistency ordering constraints, and resource 
constraints dependent on the specific modeling approach.

\noindent \textit{c. Active imaging}

The disjunctive formulation is the only one in which active imaging has been considered, that is, where assumption \hyperref[tab:constraints]{A2} is relaxed and observation duration becomes a decision variable rather than a fixed parameter
\citep{Wolfe2000,Yang2018,Lu2021,Zhou2023,Chang2025a}. This is consistent with the structural flexibility of 
continuous-time formulations, where observation duration can 
be introduced as an additional continuous variable without 
altering the overall formulation structure. 

A continuous decision variable $\text{d}^o_{ijk} \geq 0$ 
is therefore introduced to represent the observation duration 
of strip $j$ by satellite $i$ in orbit $k$. This variable is usually bounded below and above by a 
minimum $d^{\text{min}}_{ijk}$ and maximum 
$d^{\text{max}}_{ijk}$ admissible duration, respectively, 
and several definitions of these bounds have been proposed 
in the literature.
The upper bound on the observation duration is commonly defined as the maximum time the satellite can image the target within its VTW, i.e.,
\begin{align}d^{\max}_{ijk} = {vwe}_{ijk} - {vws}_{ijk},\end{align}
which ensures that any scheduled observation fits entirely within the VTW~\citep{Chang2020,Lu2021,Chang2025a}. The lower bound $d^{\min}_{ijk}$ reflects the minimum useful imaging duration. Most authors treat it as a fixed pre-defined parameter \citep{Wolfe2000,Yang2018,Lu2021,Zhou2023,Chang2025a}. \citet{Chang2020} instead derive it from the degree of competition
between targets: each target is assigned a congestion measure, and its
minimum duration decreases as this measure grows, down to a fixed
floor. A target with no conflicts is assigned a minimum duration equal
to the length of its VTW, and is therefore observed for the whole
window. \citet{Lu2022a,Lu2022b,Lu2023} define the minimum 
observation duration $d^{\min}_{ijk}$ as the minimum time 
required for the satellite to observe strip $j$ from its first 
observed endpoint to its second at maximum observation speed, which can be defined as the transition time between the two endpoints. 
In either case, new timing constraints~\eqref{eq:cont_dur_bounds} are required.
\begin{align}
    d^{\min}_{ijk} \cdot \text{x}_{ijk} \leq \text{d}^o_{ijk} 
    \leq d^{\max}_{ijk} \cdot \text{x}_{ijk}, \quad 
    \forall\ i \in \mathbf{I},\ j \in \mathbf{J},\ 
    k \in \mathbf{K}_{ij}
    \label{eq:cont_dur_bounds}
\end{align}

We note that introducing $\text{d}^o_{ijk}$ as a continuous decision variable does not affect the linearity of the linear constraints in the general model. Although the term $\text{d}^o_{ijk} \cdot \text{x}_{ijk}$ is bilinear in general, the bound constraints~\eqref{eq:cont_dur_bounds} enforce $\text{d}^o_{ijk}$ equal to 0 whenever $\text{x}_{ijk}$ equals to 0, so the product $\text{d}^o_{ijk} \cdot \text{x}_{ijk}$ reduces to $\text{d}^o_{ijk}$ throughout the feasible region, and the constraint remains linear. However, the nonlinear constraints in the general model, namely the transition time constraints~\eqref{eq:cont_nooverlap} under assumption~\hyperref[tab:constraints]{A1c}, become functions of both $\text{t}_{ijk}$ and $\text{d}^o_{ijk}$, since the transition time now depends on the end time $\text{t}_{ijk} + \text{d}^o_{ijk}$ rather than $\text{t}_{ijk} + dur^o_j$.

\noindent \textit{d. Different observation requirements}

\citet{Lemaitre2002,Cordeau2005,Habet2010} model the 
specific characteristic of free imaging directions of agile satellites, where for each target, the 
observation can be performed in any 
direction. Since after decomposition each target 
$n \in \mathbf{N}$ is associated with a fixed strip 
direction $\theta_n \in \Theta = [0^\circ, 360^\circ]$, 
each strip $j \in \mathbf{J}$ can be imaged in two 
possible acquisitions: from endpoint $(j,0)$ to endpoint 
$(j,1)$, or from endpoint $(j,1)$ to endpoint $(j,0)$.
To accommodate this, we consider the set of strip 
acquisitions $\mathbf{J}^{\text{acq}} = \{1,\ldots, 
J^{\text{acq}}\}$, with $J^{\text{acq}} = 2J$, 
associating two acquisitions with each strip 
$j \in \mathbf{J}$. By convention, acquisition $2j-1$ 
goes from endpoint $(j,0)$ to endpoint $(j,1)$ and acquisition $2j$ goes 
from endpoint $(j,1)$ to endpoint $(j,0)$. 
Constraints~\eqref{eq:cont_acq_alloc} then extend the 
allocation constraints~\eqref{eq:cont_alloc} to ensure 
that each strip is observed at most once considering both 
possible acquisitions.
\begin{align}
    & \sum_{i \in \mathbf{I}} \sum_{k \in 
    \mathbf{K}_{i(2j-1)}} \text{x}_{i(2j-1)k} + 
    \sum_{i \in \mathbf{I}} \sum_{k \in 
    \mathbf{K}_{i(2j)}} \text{x}_{i(2j)k} \leq 1, 
    \quad \forall\ j \in \mathbf{J} 
    \label{eq:cont_acq_alloc}
\end{align}
\citet{Lemaitre2002,Cordeau2005,Habet2010} additionally
consider stereoscopic requests, each requiring two
observations of the same target acquired in the same
orbital direction. To model this, let
$\mathbf{J}^{\text{ste}} \subseteq \mathbf{J}$ denote the
subset of strips belonging to a stereoscopic request. For
each $j \in \mathbf{J}^{\text{ste}}$, a twin strip $j^*$
is added to the strip set to represent the second required
observation.
Constraints~\eqref{eq:cont_stereo_alloc} then link each
pair $(j,j^*)$, enforcing that both strips are acquired in
the same orbital direction, or that neither is acquired.
\begin{subequations}
    \label{eq:cont_stereo_alloc}
\begin{align}
    &\sum_{i \in \mathbf{I}}
    \sum_{k \in \mathbf{K}_{i(2j-1)}} \text{x}_{i(2j-1)k}
    = \sum_{i \in \mathbf{I}}
    \sum_{k \in \mathbf{K}_{i(2j^*-1)}}
    \text{x}_{i(2j^*-1)k},
    \quad \forall\ j \in \mathbf{J}^{\text{ste}} \\
    &\sum_{i \in \mathbf{I}}
    \sum_{k \in \mathbf{K}_{i(2j)}} \text{x}_{i(2j)k}
    = \sum_{i \in \mathbf{I}}
    \sum_{k \in \mathbf{K}_{i(2j^*)}}
    \text{x}_{i(2j^*)k},
    \quad \forall\ j \in \mathbf{J}^{\text{ste}}
\end{align}
\end{subequations}
\citet{Xing2026} consider multiple observations of the 
same target, but unlike the already discussed approaches of 
\citet{Wang2019b,Wang2021b,Wang2025b}, each target 
requires a specific number of observations. This allows a subset of tasks to be associated 
with each target, whose cardinality equals the required 
number of observations, analogously to how 
stereoscopic requests generate a fixed set of twin 
strips. The model then remains practically 
identical to the general model up to 
notational adjustments.

A further observation requirement that has been 
considered in this class of formulations is image 
quality. Beyond the time-dependent profits based on 
image quality already discussed, additional constraints have also 
been proposed to enforce a minimum image quality level 
\citep{Liu2017,He2018,Du2020}. In this case, 
constraints~\eqref{eq:cont_image_qual} are defined to ensure that the 
image quality of the observation of strip $j$ by 
satellite $i$ in orbit $k$, determined by the 
time-dependent instant image quality function $\text{qual}(\text{t}_{ijk})$ 
evaluated at the observation start time, meets a minimum 
threshold $q_j$ (\hyperref[tab:constraints]{I5}).
\begin{align}
    &\text{qual}(\text{t}_{ijk}) \geq q_j - 
    \mathcal{M}(1 - \text{x}_{ijk}), 
    \quad \forall\ i \in \mathbf{I},\ j \in \mathbf{J},\ 
    k \in \mathbf{K}_{ij}
    \label{eq:cont_image_qual}
\end{align}
Since $\text{qual}(\cdot)$ is a nonlinear 
function of the continuous decision variable 
$\text{t}_{ijk}$, these constraints introduce 
nonlinearity into the formulation. 

\noindent \textit{e. Scheduling without orbit decomposition}

As noted in \hyperref[sec:problem_desc]{Section \ref{sec:problem_desc}}, decomposing tasks by satellite 
orbit is a modeling convention adopted in the literature to simplify 
the structure of the scheduling problem, particularly the definition 
of resource constraints. However, some authors do not consider this 
decomposition and index VTWs directly~\citep{He2018,Cho2018,Kim2026}. 
In this alternative approach, the set of VTWs of satellite~$i$ for 
strip~$j$ can be defined as $\mathbf{VW}_{ij} = \{vw_{ij1}, \ldots, vw_{ijr}, \ldots, vw_{ij(VW_{ij})}\}$, without explicit reference to the orbit in which 
each window falls. Since there is a bijection between this enumeration 
and the orbit-based indexing, where each VTW $vw_{ijr}$ corresponds to a 
unique orbit $k \in \mathbf{K}_{ij}$ such that $vw_{ijk} = vw_{ijr}$, 
the two indexing schemes are interchangeable, and the resulting 
formulation is structurally identical to the one presented above, 
with index~$k$ replaced by~$r$ throughout. 
For example, if considering the scheduling problem for passive imaging agile satellites, meaning free imaging directions, as discussed above, this formulation would not require new allocation constraints, as each strip would already be associated with two VTWs per orbit corresponding to the two possible acquisitions, and the existing allocation constraints would ensure that at most one of all the VTWs is selected for each strip. 
When resource constraints are considered, however, the two approaches are no longer equivalent. The orbit decomposition allows energy and memory budgets to be enforced independently per orbit, reflecting the physical assumption that resources are replenished at the start of each orbit. Without this decomposition, resources must instead be tracked sequentially across the entire scheduling horizon, requiring auxiliary variables to record the resource level at each event point and constraints to ensure that cumulative consumption never exceeds the available capacity at any point in time. We refer to \citet{Cho2018,Kim2026} for more details on this approach, as it is less common in the literature and falls outside the scope of this review.

\subsubsection{Arc flow}
\label{sec:arc_flow_continuous}

The main references for this formulation are \citep{Li2024,Rocha2025}.
The continuous-time arc flow formulation shares the directed-graph 
structure of the discrete arc flow formulation of
\hyperref[sec:arc_flow_local]{Section \ref{sec:arc_flow_local}}, modeling the scheduling problem for 
each satellite and orbit as a path-selection problem on a directed 
graph in which observations are represented as nodes and a directed 
arc is selected if one observation immediately precedes the other. 
The main decision variable is $\text{x}_{ijj'k}\in\{0,1\}$, equal 
to~$1$ if strip~$j'$ is observed immediately after strip~$j$ by 
satellite~$i$ in orbit~$k$, and $0$ otherwise, and the continuous 
variable $\text{t}_{ijk}\geq 0$ represents the start time of the 
observation of strip~$j$ by satellite~$i$ in orbit~$k$. For the 
consistent definition of the arc variables, a dummy source task 
$j^{s}_{ik}$ and a dummy sink task $j^{e}_{ik}$ are introduced for 
each satellite $i\in \mathbf{I}$ and orbit $k\in \mathbf{K}_i$, , with zero duration and zero 
transition times to and from any other task, so that all transitions involving 
them are feasible by construction. 

\subsubsubsection{Allocation constraints}
\noindent Constraints~\eqref{eq:ct_af_alloc} ensure that each strip is
observed at most once (\hyperref[tab:constraints]{I2}). Constraints~\eqref{eq:ct_af_source_sink} and ~\eqref{eq:ct_af_flow} are the flow conservation constraints, where constraints \eqref{eq:ct_af_source_sink}
ensure that the source and
sink tasks are always selected, and constraints~\eqref{eq:ct_af_flow} ensure that each selected observation has
exactly one immediate successor and predecessor (\hyperref[tab:constraints]{I1}).
\begin{align}
  &\sum_{i\in \mathbf{I}}\;\sum_{k\in \mathbf{K}_{ij}}\;
  \sum_{{j'\in \mathbf{J}\cup\{j^{e}_{ik}\} \setminus\{j\}}}
  \text{x}_{ijj'k}
  \;\leq\; 1,
  \quad \forall\ \,j\in \mathbf{J}
  \label{eq:ct_af_alloc}\\
  &\sum_{j\in \mathbf{J}\cup\{j^{e}_{ik}\}} \text{x}_{ij^{s}_{ik}jk}
  \;=\;
  \sum_{j\in \mathbf{J}\cup\{j^{s}_{ik}\}} \text{x}_{ijj^{e}_{ik}k}
  \;=\; 1,
  \quad \forall\ \ i\in \mathbf{I},\; k\in \mathbf{K}_i
  \label{eq:ct_af_source_sink}
  \\[4pt]
  &\sum_{{j'\in \mathbf{J}\cup\{j^{e}_{ik}\}\setminus\{j\}}} \text{x}_{ijj'k}
  \;=\;
  \sum_{{j'\in \mathbf{J}\cup\{j^{s}_{ik}\}\setminus\{j\}}} \text{x}_{ij'jk},
  \quad \forall\ \ i\in \mathbf{I},\; k\in \mathbf{K}_i,\; j\in \mathbf{J}
  \label{eq:ct_af_flow}
\end{align}
 
\subsubsubsection{Timing constraints}

\noindent Constraints~(\ref{eq:ct_af_vtw_lb}) and~(\ref{eq:ct_af_vtw_ub}) ensure that each observation 
is scheduled within its corresponding VTW (\hyperref[tab:constraints]{I3}) \citep{Li2024}.
\begin{align}
    &\text{t}_{ijk} \geq {vws}_{ijk} \cdot 
    \sum_{j'\in \mathbf{J}\cup\{j^{e}_{ik}\} \setminus\{j\}} 
    \text{x}_{ijj'k}, \quad \forall\, i\in \mathbf{I},\; j\in \mathbf{J},\; k\in \mathbf{K}_{i}
    \label{eq:ct_af_vtw_lb} \\[6pt]
    &\text{t}_{ijk} \leq \left({vwe}_{ijk} - {dur}^{{o}}_{j}\right)
    \cdot \sum_{j'\in \mathbf{J}\cup\{j^{e}_{ik}\} \setminus\{j\}} 
    \text{x}_{ijj'k}, \quad \forall\, i\in \mathbf{I},\; j\in \mathbf{J},\; k\in \mathbf{K}_{i}
    \label{eq:ct_af_vtw_ub}
\end{align}
We note that constraints~\eqref{eq:ct_af_vtw_lb} and~\eqref{eq:ct_af_vtw_ub} 
could equivalently be written with a big-${M}$ parameter, analogously 
to constraints~\eqref{eq:cont_lb} and~\eqref{eq:cont_ub} in the disjuntive 
formulation; we follow \citet{Li2024} and adopt the tighter formulation.
We also note that if the domain of $\text{t}_{ijk}$ is restricted by construction to 
$[{vws}_{ijk},\, {vwe}_{ijk} - {dur}^{{o}}_{j}]$, as considered by \citet{Rocha2025}, then the timing constraints are not required.

\subsubsubsection{Sequencing constraints}
\noindent Constraints~\eqref{eq:ct_af_trans} ensure that between two consecutively 
observed strips the corresponding transition time is respected (\hyperref[tab:constraints]{I4}). The 
following are stated for the case of sequence and time-dependent 
transition times \citep{Rocha2025}. When transition times are assumed to be sequence-dependent only (\hyperref[tab:constraints]{A1b}), 
$\text{trans}(\text{t}_{ijk},\,\text{t}_{ij'k})$ 
reduces to the fixed parameter $\text{trans}^{{o2o}}_{ijj'k}$ and 
constraints~\eqref{eq:ct_af_trans} become linear \citep{Li2024}.
\begin{align}
    &\text{t}_{ij'k} \geq\; \text{t}_{ijk} + {dur}^{{o}}_{j} 
    + \text{trans}(\text{t}_{ijk}+ dur^o_{j},\,\text{t}_{ij'k}) 
    - \mathcal{M}\left(1 - \text{x}_{ijj'k}\right), 
    \label{eq:ct_af_trans}\\
    &\quad \forall\ \ i\in \mathbf{I},\;
    j\in \mathbf{J}\cup\{j^{s}_{ik}\},\;
    j'\in \mathbf{J}\cup\{j^{e}_{ik}\}\setminus\{j\},\;
    k\in \mathbf{K}_{i} \nonumber
\end{align}
 
\subsubsubsection{Resource constraints}

\noindent Constraints~\eqref{eq:ct_af_energy} ensure that the total energy
consumed by observations and attitude maneuvers in each orbit does
not exceed the available energy budget~(\hyperref[tab:constraints]{R1}), and constraints~\eqref{eq:ct_af_memory} ensure that the total memory
consumed by observations in each orbit does not exceed the available
memory capacity~(\hyperref[tab:constraints]{R2}) \citep{Li2024}.
\begin{align}
    &\sum_{j \in \mathbf{J}} \sum_{j'\in \mathbf{J}\cup\{j^{e}_{ik}\}\setminus\{j\}} \bar{E}_i^{o} \cdot dur_{j}^o 
    \cdot \text{x}_{ijj'k} + \nonumber \\
    &   \ \sum_{j \in \mathbf{J}\cup\{j^{s}_{ik}\}} \sum_{j'\in \mathbf{J}\cup\{j^{e}_{ik}\}\setminus\{j\}}
    \bar{E}_i^{a} \cdot \text{trans}(\text{t}_{ijk} + dur^o_{j},\text{t}_{ij'k}) \cdot \text{x}_{ijj'k}
    \leq E^{\max}_i, \quad \forall \ \ i \in \mathbf{I}, k \in \mathbf{K}_i \label{eq:ct_af_energy}\\
    &\sum_{j \in \mathbf{J}\cup\{j^{s}_{ik}\}} \sum_{j'\in \mathbf{J}\cup\{j^{e}_{ik}\}\setminus\{j\}}
  \bar{M}^{o}_{i}\cdot dur^{o}_{j}\cdot \text{x}_{ijj'k}
  \;\leq\; M^{\max}_{i}, \quad \forall \ \ i \in \mathbf{I}, k \in \mathbf{K}_i
  \label{eq:ct_af_memory}
\end{align}

\subsubsubsection{Extensions}

\noindent\textit{a. Objective function variations}

\citet{Li2024} considers that each 
target requires exactly two observations. Since the usefulness of the second observation depends on how much
time has passed since the first, the authors define profit as a
function of the number of orbits elapsed between the two observations.
As this interval is itself an outcome of the scheduling decisions, an
auxiliary binary variable is introduced to identify it and select the
corresponding profit value, together with constraints that determine
the interval and ensure that profit is collected only when both
observations of a target are scheduled.
\citet{Rocha2025} considers the 
integrated scheduling problem and extend the objective function to 
jointly maximize the profit collected from observed targets and 
the total amount of data successfully downloaded to ground 
stations. The second component is modeled as the negative of 
the data volume remaining on board at the end of the planning 
horizon, penalizing observations whose data was never 
transmitted to a ground station.

\noindent\textit{b. Integrated scheduling problem.}

\citet{Rocha2025} extend the arc flow formulation to the integrated
observation and downlink scheduling problem. Differently from the
formulations described above, the authors index all tasks directly by
their time window, without reference to an explicit orbit structure.
Observation and downlink tasks are treated as nodes of a single
directed graph per satellite, and a unified binary arc variable is
defined for every pair of consecutive tasks, irrespective of their
type. Flow conservation constraints are imposed at every node, so that
each satellite follows a single sequence of tasks over the whole
planning horizon.
The authors also adopt a different coupling between observation and
downlink decisions. Downlink variables carry no strip index, and
onboard storage is instead modeled as a level that increases during
observation tasks and decreases during downlink tasks, whose duration
is itself a decision variable. The incentive to transmit is provided
by penalizing the data volume remaining onboard at the end of the
planning horizon. This removes constraints
\hyperref[tab:constraints]{ID1} and \hyperref[tab:constraints]{ID2}
and reduces the number of binary variables, but the model then
guarantees a total transmitted volume rather than the delivery of
individual requests, so that requirements defined per request, such as
a deadline for the data of a specific target, cannot be expressed.
Since this formulation does not adopt the orbit-decomposition structure on which our paper is based and its modeling details differ substantially from those described so far, we refer to \citet{Rocha2025} for a more detailed description.

Finally, we note that this formulation has not yet been extended 
to model active imaging satellites, where observation duration 
is treated as a decision variable.

\section{Conclusions and future directions}
\label{sec:conclusions}

This paper has presented a comprehensive review of optimization models
for EOS scheduling. A classification of EOS scheduling problems was proposed, identifying the constraints and assumptions that define them and the problem features that have a direct consequence on the structure of the resulting model, and expressing each problem through a three-field notation.
A classification of optimization
models, based on their time and event representations, was then used as
a framework for describing the main formulations proposed in the
literature. Modeling aspects of representative formulations of each
category were presented in a unified notation, emphasizing their main
variables and constraints, the extensions proposed for each, and their
strengths and limitations. Finally, the structural trade-offs between
the alternative representations were discussed, together with the
modeling issues that remain to be solved across all of them.

The first classification concerns the problem and identifies the
features that have a direct consequence on the structure of the
resulting model, namely mission scope, satellite agility, transition
times, observation requirements, resource constraints and objective
function. Within this classification we proposed a new taxonomy for
satellite agility, grounded in two physical attributes, axis mobility
and attitude control during observations, which replaces the
conventional distinction between conventional, agile and super-agile
satellites.
In answer to RQ1, what
distinguishes one EOS scheduling problem from another at the level of
the model is therefore a small set of features, of which the decisive
ones are those that determine whether the timing of a task, namely the start time and duration, is a
decision, since it is that property that governs whether the
parameters of the model can be pre-computed and the formulation
remains linear. In the proposed notation, these features correspond mainly to the satellite type $\alpha_1$, which determines whether start times and durations are decisions, and to the treatment of transition times $\beta_2$.

The second classification concerns the models and is based on their
time representation, discrete or continuous, and on their event
representation, which determines how the arrangement of tasks over
time is encoded. Reviewing the formulations along these two dimensions
shows that each representation handles the time dependence of the
problem at a different cost, either in the mathematical structure of
the constraints or in the size of the model. It also shows that the
extensions proposed in the literature are unevenly distributed across
the categories, with active imaging modeled only in the disjunctive
continuous-time formulation, the integrated imaging and downlink
problem addressed in every category except those based on local time
discretization, and image quality constraints appearing almost
exclusively in continuous-time models. In answer to RQ2, no
representation dominates the others. The time representation
determines whether the parameters that depend on task timing can be
pre-computed, and therefore whether the model remains linear, while
the event representation determines not only the size of the model but
which problem features can be more easily expressed, which is what explains
the distribution of the extensions.

Several directions for future research follow directly from this
analysis.
The most immediate one concerns the exact treatment of sequence and
time-dependent transition times, which no reviewed formulation
achieves at the resolution at which the underlying orbital data are
available. Existing approaches either discretize the visibility
windows, which limits the size of the instances that can be solved, or
approximate the angle profiles and the resulting transition times,
which makes the model tractable at the cost of an approximation whose
effect on solution quality is not quantified. The structural
properties of the transition time function, in particular the
first-in-first-out rule and the triangle inequality presented by
\citet{Peng2022}, have so far been exploited only within heuristic
search, and using them to derive exact formulations appears to be a
promising avenue.
A second direction concerns the
integrated scheduling problem of active imaging agile satellites. Since the observation duration
determines the volume of data generated by each observation, treating
it as a decision variable couples the imaging and the downlink
decisions in a way that does not arise when durations are fixed, and
this coupling extends to the completion times on which the causality
and the transition constraints of the integrated problem depend. No existing formulation captures it. Closely related to this, the
representation of onboard resources deserves more attention, as energy
consumed by maneuvers, energy recharged by solar panels and memory
carried across successive orbits are usually treated through
simplifications that are convenient for the model but that have not
been assessed against more detailed alternatives.
Finally, a third direction concerns the empirical comparison of the
formulations. This paper establishes how the alternative
representations differ structurally, but their relative computational
performance remains unknown, since studies proposing a formulation
compare it against the algorithms of the works they extend rather than
against alternative formulations of the same problem, and the
instances used differ in size, geographic distribution and modeling
assumptions. Public benchmark instances reporting the orbital data and
the associated look angle profiles, which are required to evaluate
time-dependent transition times, would make such a comparison possible
and would allow the structural trade-offs discussed here to be
quantified.

In sum, this review provides a unified framework for understanding the EOS scheduling literature by harmonizing problem definitions, modeling assumptions, and optimization formulations within a consistent taxonomy. By clarifying the structural relationships between problem characteristics and model design, it establishes a common basis for comparing existing formulations and interpreting their trade-offs. At the same time, by identifying the main unresolved modeling challenges, it offers a clearer foundation for the development and evaluation of future optimization models in EOS scheduling.

\section*{CRediT authorship contribution statement}
\textbf{Margarida Caleiras:} Conceptualization, Methodology, Investigation, Formal analysis, Visualization, Writing -- original draft.
\textbf{Paulo Jorge Nascimento:} Conceptualization, Methodology, Supervision, Writing -- review \& editing.
\textbf{Samuel Moniz:} Conceptualization, Methodology, Supervision, Funding acquisition, Writing -- review \& editing.

\section*{Data availability}
No data was used for the research described in the article.

\section*{Declaration of generative AI and AI-assisted technologies in the manuscript preparation process}
During the preparation of this work, the authors used Claude in
order to improve language and readability. After using this tool/service,
the authors reviewed and edited the content as needed and took full
responsibility for the content of the publication.

\section*{Acknowledgements}
This research is sponsored by national funds through FCT – Fundação para a Ciência e a Tecnologia, under projects UID/00285/2025 and LA/P/0112/2020, and by the European Union under the Next Generation EU, through a grant of the Portuguese Republic's Recovery and Resilience Plan (RRP) Partnership Agreement, within the scope of the project PRODUTECH R3 – "Agenda Mobilizadora da Fileira das Tecnologias de Produção para a Reindustrialização".
All supports are gratefully acknowledged. The funding sources had no involvement in the study design, the analysis and interpretation of the literature, the writing of the article, or the decision to submit it for publication.

\printbibliography
\end{document}